\documentclass[12pt,a4paper]{article}

\usepackage[english]{babel}
\usepackage[dvips]{graphicx}

\usepackage{yhmath}
\usepackage{amsfonts}
\usepackage{amsmath}
\usepackage[single]{accents}

\usepackage{array}
\usepackage{verbatim}
\usepackage{graphics}
\usepackage{subfigure}
\usepackage{epsfig}
\usepackage{tabularx}
\usepackage{mathcomp}
\usepackage{url}

\newcommand{\mB}{\mathcal{B}}

\newcommand{\mG}{\mathcal{G}}

\newcommand{\mS}{\mathcal{S}}
\newcommand{\mT}{\mathcal{T}}

\newcommand{\mV}{\mathcal{V}}

\newcommand{\Lij}{\mathcal{L}_{ij}}

\newcommand{\Li}{\mathcal{L}_{i}}

\newcommand{\Bijrs}{\mB_{ijrs}}
\newcommand{\Bir}{\mB_{ir}}

\newcommand{\Cab}{\mathcal{C}_{\alpha\beta}}
\newcommand{\Ca}{\mathcal{C}_{\ga}}

\newcommand{\Sij}{\mathcal{S}_{ij}}
\newcommand{\Skk}{\mathcal{S}_{kk}}

\newcommand{\bu}{\mathbf{u}}

\newcommand{\bF}{\mathbf{F}}

\newcommand{\bmG}{\boldsymbol{\mG}}

\newcommand{\bS}{\mathbf{S}}

\newcommand{\bU}{\mathbf{U}}

\newcommand{\bq}{\mathbf{q}}

\newcommand{\bbar}{\overline}

\newcommand{\frho}{\bbar{\rho}}

\newcommand{\fp}{\bbar{p}}

\newcommand{\wt}{\widetilde}

\newcommand{\fu}{\wt{u}}
\newcommand{\fh}{\wt{h}}
\newcommand{\fe}{\wt{e}}

\newcommand{\fT}{\wt{T}}

\newcommand{\fS}{\wt{\mathcal{S}}}
\newcommand{\fSij}{\wt{\mS}_{ij}}
\newcommand{\fSrs}{\wt{\mS}_{rs}}
\newcommand{\fSkk}{\wt{\mS}_{kk}}

\newcommand{\fsigmaij}{\wt{\sigma}_{ij}}

\newcommand{\wh}{\widehat}

\newcommand{\de}{\partial}
\newcommand{\Div}{\nabla\cdot}

\newcommand{\ga}{\alpha}
\newcommand{\gb}{\beta}

\newcommand{\w}{\omega}

\newcommand{\sigmaij}{\sigma_{ij}}
\newcommand{\taukk}{\tau_{kk}}
\newcommand{\tauij}{\tau_{ij}}
\newcommand{\ei}{e_{\rm i}}

\newcommand{\hfu}{\breve{\fu}}
\newcommand{\hfT}{\breve{\fT}}
\newcommand{\hfrho}{\widehat{\frho}}
\newcommand{\hfS}{\breve{\fS}}

\newcommand{\bn}{\mathbf{n}}
\newcommand{\bx}{\mathbf{x}}
\newcommand{\br}{\mathbf{r}}

\newcommand{\hdelta}{\widehat{\Delta}}

\newcommand{\deltaij}{\delta_{ij}}
\newcommand{\nusgs}{\nu^{{\rm sgs}}}

\newcommand{\RE}{Re}
\newcommand{\MA}{M\hspace{-1pt}a}
\newcommand{\PR}{Pr}
\newcommand{\dm}{{\rm d}}

\begin{document}

\title{Anisotropic dynamic models for Large Eddy Simulation of compressible flows with a high order DG method}

\author{Antonella Abb\`a $^{(1)}$, Luca Bonaventura$^ {(2)}$\\
Michele Nini$^{(1)}$, Marco Restelli$^{(3)}$}

\maketitle

 \begin{center}
{\small
(1) Dipartimento di Ingegneria Aerospaziale, Politecnico di Milano \\
Via La Masa 34, 20156 Milano, Italy\\
{\tt antonella.abba@polimi.it, michele.nini@polimi.it}\\
{$ \ \ $ }\\
(2) MOX -- Modelling and Scientific Computing, \\
Dipartimento di Matematica, Politecnico di Milano \\
Via Bonardi 9, 20133 Milano, Italy\\
{\tt luca.bonaventura@polimi.it}\\
{$ \ \ $ }\\
(3) NMPP -- Numerische Methoden in der Plasmaphysik \\
Max--Planck--Institut f\"ur Plasmaphysik \\
Boltzmannstra\ss e 2, D-85748 Garching, Germany \\
{\tt  marco.restelli@ipp.mpg.de}\\
}
\end{center}

\date{}


\vspace*{0.5cm}


\vspace*{0.5cm}

\pagebreak

\abstract{The impact of anisotropic dynamic models for  applications to LES
of compressible flows is assessed in the framework of a 
numerical model based on high order discontinuous finite elements.
 The projections onto lower dimensional subspaces associated to 
 lower degree basis function are used as LES filter, along the lines  proposed
in  Variational Multiscale templates. 
  Comparisons with DNS results available in the literature for channel flows at Mach numbers 0.2, 0.7 and 1.5
show clearly that the anisotropic model is able to reproduce well some key features of the flow,
especially close to the wall, where the flow anisotropy plays a major role.}


\pagebreak

\section{Introduction}
\label{intro} \indent

High order finite element methods are an extremely appealing framework to implement LES models of turbulent flows,
due to their potential for reducing the impact of numerical dissipation on most of the spatial scales of interest.
Discontinuous Galerkin (DG) methods have been applied to LES  and DNS by several authors,
see e.g. \cite{collis:2002b}, \cite{collis:2002a}, 
\cite{vanderbos:2007},  \cite{landmann:2008}, \cite{sengupta:2007}, \cite{uranga:2011}, \cite{vanderbos:2010},  \cite{wei:2011}. 
DG methods are particularly appealing for realistic CFD applications 
for a number of practical and conceptual reasons.  At a more practical level,
they allow to implement $h$ and $p$ refinement procedures  with great ease and to work on complex and also non conforming
meshes. Even though they imply quite stringent stability restriction for explicit time discretization approaches,
a number of techniques is available to improve computational efficiency if required, see e.g. 
\cite{garcia:2014}, \cite{giraldo:2010}, \cite{restelli:2009}, \cite{schulze:2009}, \cite{tumolo:2013}. At a more conceptual level,
discontinuous finite elements provide a natural framework to generalize LES filters to arbitrary
computational meshes.  As proposed in some of the previously quoted papers, the filter operator
that is the key tool in LES can be identified with the projection operator
on a finite dimensional space related to the discretization. This allows to generalize easily the LES
 concept to unstructured meshes and complex geometries. Ideas of this kind have first arisen in the framework
 of the so called   Variational Multiscale (VMS) approach, that has been introduced in \cite{hughes:1998} and applied to
Large Eddy Simulation (LES) of incompressible flows in \cite{hughes:2000},\cite{hughes:2001a},\cite{hughes:2001b}
(see also the review in \cite{hughes:2004}).
Other multiscale  approaches to LES in the framework of finite element discretizations
have been proposed e.g. in \cite{john:2010}, \cite{john:2007}, \cite{koobus:2004}, \cite{munts:2007}. 
 
This very promising framework, however, seems to have been only partially exploited so far.
In  \cite{vanderbos:2010}, for example, the LES filter has been realized
by face based projection operators  that are different from those for which
the VMS template has been outlined in \cite{vanderbos:2007}. 
Furthermore, to the best of our knowledge, only simple Smagorinsky closures have  
been employed  to model the subgrid stresses in the available literature.
In this paper, we investigate the potential benefit  resulting from the use of
the anisotropic dynamic model  \cite{abba:2003}, appropriately extended to the compressible case,
 in the context of a high order DG numerical model. Anisotropic models try to address the 
failure of the Boussinesq hypothesis (see e.g. \cite{schmitt:2007} for an extensive review of this subject)
by introducing a tensor valued subgrid viscosity, thus avoiding alignment of the stress and velocity strain rate tensors.
  We implement a LES model with projection-based filter
   in the framework of a high order DG  method
and we assess the performance of this more sophisticated subgrid closure  with respect to the simple Smagorinsky 
closure. The comparison is carried out with respect to the DNS
experiment results reported in  \cite{coleman:1995}, \cite{moser:1999} and \cite{wei:2011}. 
The results of the comparison show a clear improvement in the prediction of several key features of the flow with respect
to the Smagorinsky closure implemented in the same framework. In particular, the anisotropic model
allows to achieve a better representation of mean  profiles, turbulent stresses and, more generally, 
of the total turbulent kinetic energy. 
The proposed approach appears to lead to significant improvements also
in the lower Mach number regimes, which justifies further extensions to flows in presence of gravity, 
with the goal of providing turbulence 
models for applications to environmental stratified flows that do not 
require {\it ad hoc} tuning of parameters. Furthermore, the numerical framework that is validated in this
paper will be employed for the assessment of the proposal presented in \cite{germano:2014} for the
extension of the eddy viscosity model to compressible flows.

In section \ref{modeq}, the Navier-Stokes equations for compressible flow are recalled. 
In section \ref{subgrid}, the LES models employed are described.
In section \ref{dg}, the DG finite element discretization is reviewed,
while in section \ref{results} the results of our comparisons with DNS data are reported.
Some conclusions and perspectives for future work are presented in section \ref{conclu}.

\section{Model equations}
\label{modeq} \indent

We consider the compressible Navier--Stokes equations,
which, employing the Einstein notation, can be written in dimensional
form (denoted by the superscript ``d'') as
\begin{subequations}
\label{eq:nscompr-dim}
\begin{align}
&\de_{t^\dm} \rho^\dm + \de_{x^\dm_j} (\rho^\dm u^\dm_j) = 0 \\
&\de_{t^\dm} (\rho^\dm u^\dm_i) + \de_{x^\dm_j} (\rho^\dm u^\dm_i
 u^\dm_j) + \de_{x^\dm_i} p^\dm - \de_{x^\dm_j} \sigmaij^\dm =
 \rho^\dm f^\dm_i \\
&\de_{t^\dm} (\rho^\dm e^\dm) + \de_{x^\dm_j} (\rho^\dm h^\dm u^\dm_j)
- \de_{x^\dm_j} (u^\dm_i \sigmaij^\dm) + \de_{x^\dm_j} q^\dm_j =
\rho^\dm f_j^\dm u_j^\dm ,
\end{align}
\end{subequations}
where $\rho^\dm$, $\bu^\dm$ and $e^\dm$ denote density, velocity and
specific total energy, respectively, $p^\dm$ is the pressure,
$\mathbf{f}^\dm$ is a prescribed forcing, $h^\dm$ is the specific
enthalpy, defined by $\rho^\dm h^\dm=\rho^\dm e^\dm+p^\dm$, and
$\sigma^\dm$ and $\bq^\dm$ are the diffusive momentum and heat fluxes.
Equation~(\ref{eq:nscompr-dim}) must be complemented with the equation of state
\begin{equation}
p^\dm = \rho^\dm R T^\dm,
\label{eq:state-eq-dim}
\end{equation}
where $T^\dm$ is the temperature and $R$ is the ideal gas constant.
The temperature can then be expressed in terms of the prognostic
variables introducing the specific internal energy $\ei^\dm$, so that
\begin{equation}
e^\dm = \ei^\dm + \frac{1}{2}u^\dm_ku^\dm_k,
\qquad \ei^\dm = c_v T^\dm,
\label{eq:ei-dim}
\end{equation}
where $c_v$ is the specific heat at constant volume. Finally, the
model is closed with the constitutive equations for the diffusive
fluxes
\begin{equation}
\sigmaij^\dm = \mu^\dm \Sij^{d,\dm}, \qquad
q^\dm_i = -\frac{\mu^\dm c_p}{\PR} \de_{x^\dm_i} T^\dm,
\label{eq:constitutive-dim}
\end{equation}
where $\Sij^\dm = \de_{x^\dm_j} u_i^\dm + \de_{x^\dm_i}
u^\dm_j$ and $\Sij^{d,\dm} = \Sij^\dm -
\dfrac{1}{3}\Skk^\dm\delta_{ij}$, the specific heat at constant
pressure is $c_p=R+c_v$, $\PR$ denotes the Prandtl number,
and the dynamic viscosity $\mu^d$ is assumed to depend only on
temperature $T^d$ according to the power law
\begin{equation}\label{eqn:mu_def-dim}
\mu^\dm(T^\dm) = \mu^\dm_0 \left(\dfrac{T^\dm}{T^\dm_0}\right)^{\alpha},
\end{equation}
in agreement with Sutherland's hypothesis (see e.g.
\cite{schlichting:1979}) with $\alpha=0.7$.
The dimensionless form of the problem is obtained assuming reference
quantities $\rho_r$, $L_r$, $V_r$ and $T_r$, as well as
\begin{equation}
\begin{array}{llll}
t_r = \frac{L_r}{V_r}, & p_r = \rho_r R T_r, &
\sigma_r = \frac{\mu_rV_r}{L_r}, & f_r = \frac{V_r^2}{L_r}, \\[2mm]
e_r = R T_r, & q_r = \frac{\mu_rc_p T_r}{\PR\, L_r}, &
\multicolumn{2}{l}{
\mu_r = \mu_0^\dm\left( \frac{T_r}{T^\dm_0} \right)^\alpha.
}
\end{array}
\label{eq:reference-quantities}
\end{equation}
Defining now
\begin{equation}
\begin{array}{lll}
\rho^\dm = \rho_r\rho, & u_i^\dm = V_ru_i, & T^\dm = T_rT, \\[2mm]
t_r \de_{t^\dm} = \de_t, & L_r \de_{x^\dm_i} = \de_i, \\[2mm]
p^\dm = p_rp, & \sigmaij^\dm = \sigma_r\sigmaij, &
f^\dm = f_rf, \\[2mm]
e^\dm = e_re, & q^\dm = q_rq, & \ei^\dm = e_r\ei, \\[2mm]
\mu^\dm = \mu_r\mu,
\end{array}
\label{eq:normalized-quantities}
\end{equation}
we obtain
\begin{subequations}
\label{eq:nscompr}
\begin{align}
&\de_t \rho + \de_j (\rho u_j) = 0 \\
&\de_t (\rho u_i) + \de_j (\rho u_i u_j) + 
\frac{1}{\gamma\,\MA^2}\de_i p - \frac{1}{\RE}\de_j \sigmaij =
\rho f_i \label{eq:nscompr-momentum} \\
&\de_t (\rho e) + \de_j (\rho h u_j)
- \frac{\gamma\,\MA^2}{\RE} \de_j (u_i \sigmaij)  \nonumber \\
&\hspace{3cm}+ \frac{1}{\kappa\RE\PR}\de_j q_j = \gamma\MA^2\rho f_j u_j ,
\end{align}
\end{subequations}
where
\begin{equation}
\MA = \frac{V_r}{\left( \gamma R T_r \right)^{1/2}}, \qquad
\RE = \frac{\rho_rV_rL_r}{\mu_r}
\label{eq:adim-numbers}
\end{equation}
and
\[
\qquad \rho h=\rho e+p, \qquad
\gamma = \frac{c_p}{c_v}, \qquad
\kappa = \frac{R}{c_p}.
\]
Other relevant equations in dimensionless form are the equation of
state
\begin{equation}
p = \rho T,
\label{eq:state-eq}
\end{equation}
the definition of the internal energy 
\begin{equation}
e = \ei + \frac{\gamma\MA^2}{2} u_ku_k,
\qquad \ei = \frac{1-\kappa}{\kappa} T,
\label{eq:ei}
\end{equation}
the constitutive equations 
\begin{equation}
\sigmaij = \mu \Sij^d, \qquad
q_i = -\mu\de_i T,
\label{eq:constitutive}
\end{equation}
with $\Sij = \de_j u_i + \de_i u_j$ and $\Sij^d = \Sij -
\dfrac{1}{3}\Skk\delta_{ij}$, and the temperature dependent
viscosity 
\begin{equation}\label{eqn:mu_def}
\mu(T) = T^\alpha.
\end{equation}
In order to derive the filtered equations for the LES model, an
appropriate filter has to be introduced, which will be denoted by the
operator $\bbar{\cdot}$ and which is assumed to be characterized by a
spatial scale $\Delta$.  Using an approach that recalls the VMS
concept, the precise definition of this operator, as well as of the
associated scale, will be built in the numerical DG discretization.
Such a definition will be given in section \ref{dg}; here, we mention
that $\Delta$ will in general depend on the local element size
and therefore has to be interpreted as a piecewise constant function
in space.
As customary in compressible LES, see e.g. \cite{garnier:2009}, in order to avoid subgrid terms arising in the
continuity equation, we also introduce the Favre filtering
operator $\wt{\cdot}$, defined implicitly by the Favre decomposition
\begin{equation}\label{eqn:favre_decomp}
 \bbar{\rho u_i} = \frho \wt{u}_i, \qquad
 \bbar{\rho e} = \frho \wt{e}.
\end{equation}
Similar decompositions are introduced for the internal energy and the
enthalpy
\[
 \bbar{\rho \ei} = \frho \wt{\ei}, \qquad
 \bbar{\rho h} = \frho \wt{h} = \frho \wt{e} + \bbar{p},
\]
as well as for the temperature, which, taking into
account~(\ref{eq:state-eq}), yields
\begin{equation}
 \bbar{\rho T} = \frho \wt{T} = \bbar{p}.
\label{eq:state-eq-Favre}
\end{equation}
Equation~(\ref{eq:ei}) then  implies
\begin{equation}
 \frho \wt{e} = \frho \wt{\ei} +
 \frac{\gamma\MA^2}{2} \left( \frho\wt{u}_k\wt{u}_k 
 + \tau_{kk}\right), \qquad
\frho\wt{\ei} = \frac{1-\kappa}{\kappa} \frho\wt{T},
\label{eq:ei-Favre}
\end{equation}
where, as customary,
\begin{equation}
  \tauij = \bbar{\rho u_i u_j} - \frho\fu_i\fu_j.
\label{eqn:tauij_sgs}
\end{equation}
Notice that, from~(\ref{eq:ei-Favre}), $\tau_{kk}$ represents the
filtered turbulent kinetic energy. Finally, neglecing the subgrid scale contributions,
let us introduce a filtered counterpart of~(\ref{eq:constitutive}), namely
\begin{equation}
\fsigmaij = \mu(\wt{T}) \fSij^d, \qquad
\wt{q}_i = -\mu(\wt{T}) \de_i \fT,
\label{eq:constitutive-Favre}
\end{equation} 
with $\fSij = \de_j \fu_i + \de_i \fu_j$ and $\fSij^d = \fSij -
\dfrac{1}{3}\fSkk\delta_{ij}$. With these definitions, and disregarding
the commutation error of the filter operator with respect to space and time
differentiation, the filtered form of~(\ref{eq:nscompr}) is 
\begin{subequations}
\label{eq:nscompr-filtered-intermediate}
\begin{align}
&\de_t \frho + \de_j (\frho \fu_j) = 0 \\
&\de_t \left( \frho \fu_i \right) + \de_j \left(\frho \fu_i \fu_j\right) 
+ \frac{1}{\gamma\,\MA^2}\de_i \fp -\frac{1}{\RE} \de_j \fsigmaij 
\nonumber \\
& \qquad \qquad
= - \de_j \tauij - \de_j \epsilon^{{\rm sgs}}_{ij} + \frho f_i\\
& \de_t \left(\frho\fe\right) + \de_j \left(\frho\fh \fu_j\right) 
- \frac{\gamma\,\MA^2}{\RE}\de_j \left(\fu_i \fsigmaij \right)
+ \frac{1}{\kappa\RE\PR}\de_j \wt{q}_j \nonumber \\
& \qquad \qquad =
- \de_j \left(\rho h u_j\right)^{{\rm sgs}}
+ \frac{\gamma\,\MA^2}{\RE}\de_j \phi^{{\rm sgs}}_j \nonumber \\
& \qquad \qquad - \frac{1}{\kappa\RE\PR}\de_j \theta^{{\rm sgs}}_j + \gamma\MA^2\frho f_j \fu_j,
\end{align}
\end{subequations}
where
\begin{equation}
\begin{array}{ll}
\epsilon^{{\rm sgs}}_{ij} = \bbar{\sigma}_{ij} - \fsigmaij, &
 \left(\rho h u_i\right)^{{\rm sgs}} = \bbar{\rho hu_i} -
 \frho\fh\fu_i,
\\[2mm]
 \phi^{{\rm sgs}}_j = \bbar{u_i \sigmaij} - \fu_i\fsigmaij, &
 \theta^{{\rm sgs}}_i = \bbar{q}_i - \wt{q}_i.
\end{array}
\label{eqn:epsij_sgs}
\end{equation}
Notice that, in order to avoid unnecessary complications, and
  since this is the case for the numerical results considered in this
  work, we assume in (\ref{eq:nscompr-filtered-intermediate}) that $f_j$ is uniform in space.
Based on the analyses e.g. in \cite{piomelli:2000} and \cite{vreman:1995} and on the
fact that 
\begin{equation}
 \bbar{\sigma}_{ij} \approx \fsigmaij, \qquad
 \bbar{q}_i \approx \wt{q}_i
\end{equation}
the term $\de_j \phi^{{\rm sgs}}_j$ is considered to be negligible, as
well as $\epsilon^{{\rm sgs}}_{ij}$ and $\theta^{{\rm sgs}}_j$.
Concerning the subgrid enthalpy flux, we proceed as follows. First of
all, notice that using~(\ref{eq:state-eq}) and~(\ref{eq:ei}), as well
as their filtered counterparts~(\ref{eq:state-eq-Favre})
and~(\ref{eq:ei-Favre}), we have
\[
\rho h = \frac{1}{\kappa}\rho T + \frac{\gamma\MA^2}{2}\rho u_ku_k, \qquad
\frho\wt{h} = \frac{1}{\kappa}\frho \wt{T}
+ \frac{\gamma\MA^2}{2}\left( \frho\fu_k\fu_k + \taukk \right).
\]
Introducing now the subgrid heat and turbulent diffusion fluxes
\begin{subequations}
\begin{align}
Q_i^{{\rm sgs}} & = \bbar{\rho u_i T} - \frho\fu_i\fT =
 \frho \left( \wt{u_i T} - \fu_i\fT \right) \label{eqn:Qj_sgs} \\
J_i^{{\rm sgs}} & = \bbar{\rho u_iu_ku_k} - \frho \fu_i\fu_k\fu_k =
\frho \wt{u_i u_k u_k} - \frho\fu_i\fu_k\fu_k
 \label{eqn:Jj_sgs}
\end{align}
\end{subequations}
we have
\begin{equation}
\left(\rho h u_i\right)^{{\rm sgs}} = \frac{1}{\kappa} Q_i^{{\rm sgs}}
+ \frac{\gamma\MA^2}{2} \left( J_i^{{\rm sgs}} - \taukk\fu_i \right).
\label{eq:enthalpy-sgs}
\end{equation}
Notice that, introducing the generalized central moments
$\tau(u_i,u_j,u_k)$ as in \cite{germano:1992}, with
\begin{equation}
\tau(u_i,u_j,u_k) = \frho\wt{u_i u_j u_k} - \fu_i\tau_{jk} 
  - \fu_j\tau_{ik} - \fu_k\tauij - \frho\fu_i\fu_j\fu_k,
\end{equation}
$J_i^{{\rm sgs}}$ in~(\ref{eqn:Jj_sgs}) can be rewritten as
\begin{equation}\label{eqn:Jj_sgs2}
 J_i^{{\rm sgs}} = \tau(u_i,u_k,u_k) + 2\fu_k\tau_{ik} + \fu_i\taukk.
\end{equation}
Summarizing, given the above approximations and definitions, the
filtered equations~(\ref{eq:nscompr-filtered-intermediate}) become
\begin{subequations}
\label{filteq}
\begin{align}
&\de_t \frho + \de_j (\frho \fu_j) = 0 \\
&\de_t \left( \frho \fu_i \right) + \de_j \left(\frho \fu_i \fu_j\right) 
+ \frac{1}{\gamma\,\MA^2}\de_i \fp -\frac{1}{\RE} \de_j \fsigmaij
= - \de_j \tauij + \frho f_i \label{filteq-momentum} \\
& \de_t \left(\frho\fe\right) + \de_j \left(\frho\fh \fu_j\right) 
- \frac{\gamma\,\MA^2}{\RE}\de_j \left(\fu_i \fsigmaij \right)
+ \frac{1}{\kappa\RE\PR}\de_j \wt{q}_j \nonumber \\
& \qquad \qquad =
- \frac{1}{\kappa}\de_j Q_j^{{\rm sgs}}
- \frac{\gamma\MA^2}{2}\de_j \left( J_j^{{\rm sgs}} - \taukk\fu_j \right) \nonumber \\
& \qquad \qquad+ \gamma\MA^2\frho f_j \fu_j.
\label{filteq-energy} 
\end{align}
\end{subequations}

\section{Subgrid models}
\label{subgrid} \indent
We will now introduce the subgrid models used in our LES experiments.
Firstly, we will briefly recall the formulation of the classical Smagorinsky subgrid model,
which, in spite of its limitations (see e.g. the discussion in \cite{sagaut:2006}), has been applied almost
exclusively in the DG-LES models proposed in the literature so far.
Then, we will discuss a dynamic, anisotropic subgrid model proposed
in~\cite{abba:2003} that does not suffer from various limitations of the
Smagorinsky model and  is here extended to the compressible case.

\subsection{The Smagorinsky model}
\label{subgrid:smagorinsky}

 In a Smagorinsky-type model, the deviatoric part of the subgrid stress tensor $\tau_{ij}$ in (\ref{filteq})
is modelled by an isotropic, scalar turbulent viscosity
$\nu^{{\rm sgs}}$, yielding
\begin{subequations}\label{eqn:nu_smag}
\begin{align}
 &  \tauij -\frac{1}{3}\tau_{kk} \deltaij = - \frac{1}{\RE} \frho \nu^{{\rm sgs}} \fSij^d,
\label{eqn:nu_smag:tauij}
 \\
& \nu^{{\rm sgs}} = \RE\, C_S^2\Delta^2 |\fS| f_D,
\label{eqn:nu_smag:nu}
\end{align}
\end{subequations}
where $C_S=0.1$ is the Smagorinsky constant, $|\fS|^2 =
\dfrac{1}{2}\fSij\fSij$ and $\Delta$ is the filter scale introduced in
section~\ref{modeq}.  The Van Driest damping function
in~(\ref{eqn:nu_smag:nu}) is defined as
\begin{equation}
f_D(y^+) = 1 - \exp\left(- y^+/A \right),  
\end{equation}
where $A$ is a constant and $y^+=\frac{\rho_r u_\tau^\dm d^\dm_{{\rm
wall}}}{\mu_r}$, with $d^\dm_{{\rm wall}}$ denoting the (dimensional)
distance from the wall and $u_\tau^\dm$ the (dimensional) friction
velocity. The
introduction of such a damping function in~(\ref{eqn:nu_smag:nu})  is
necessary to reduce the length $\Delta$ according to the smaller size
of turbulent structures close to the wall and to recover the correct
physical trend for the turbulent viscosity (see for instance
\cite{sagaut:2006}); in the following, the value $A=25$ is employed.
We also notice that the Reynolds number has been included in the
definition of $\nu^{{\rm sgs}}$ so that the corresponding dimensional
viscosity can be obtained as $\nu^{{\rm sgs,d}}=\frac{\mu_r}{\rho_r}
\nu^{{\rm sgs}}$.

Concerning the isotropic part of the subgrid stress tensor,
some authors \cite{erlebacher:1992} have neglected it,
considering it negligible with respect to the pressure contribution.
Alternatively, following  \cite{yoshizawa:1986}, the isotropic components of the subgrid stress tensor can be
modelled as:
\begin{equation}\label{eqn:taukk_smag}
 \taukk = C_I \frho\Delta^2 |\fS|^2.
\end{equation}
Along the lines of \cite{eidson:1985}, the subgrid temperature
flux~\eqref{eqn:Qj_sgs} is assumed to be proportional 
to the resolved temperature gradient and is modelled with the eddy viscosity model
\begin{equation}\label{eqn:Qj_smag}
 Q_i^{{\rm sgs}} = - \frac{\PR}{\PR^{{\rm sgs}}} \frho \nusgs \de_i\fT,
\end{equation}
where $\PR^{{\rm sgs}}$ is a subgrid Prandtl number. Notice that the
corresponding dimensional flux is $Q_i^{{\rm sgs,d}} = q_r Q_i^{{\rm
sgs}}$. Finally, concerning $J_i^{{\rm sgs}}$ in (\ref{eqn:Jj_sgs2}), by
analogy with RANS models, the term $\tau(u_i,u_j,u_k)$ is neglected
(see e.g. \cite{knight:1998}), yielding
\begin{equation}\label{eqn:Jj_smag}
 J_i^{{\rm sgs}} \approx 2\fu_k\tau_{ik} + \fu_i\taukk.
\end{equation}

\subsection{The anisotropic model}
\label{subgrid:anisotropic}

We consider now the dynamic, anisotropic subgrid model proposed in
\cite{abba:2003}, which is extended here to the compressible case.
This approach has the goal of removing two limitations of the
Smagorinsky model: the fact that the constants $C_S$ and $C_I$ must be
chosen \emph{a priori} for the whole domain, and the alignment of the
subgrid flux tensors with the gradients of the corresponding
quantities. The first limitation is removed employing the Germano
dynamic procedure~\cite{germano:1991}, while the subgrid tensor
alignment is removed generalizing the proportionality relations such
as~(\ref{eqn:nu_smag:tauij}) introducing proportionality parameters
which are tensors rather than scalar quantities.

More specifically,
 the subgrid stress tensor $\tauij$ is assumed proportional to the strain rate tensor through a fourth order symmetric
 tensor as follows
 \begin{equation}\label{eqn:tauij_aniso}
 \tauij =  -\frho \Delta^2 |\fS| \Bijrs\fSrs.
\end{equation}
To compute dynamically the tensor $\Bijrs$, let us first observe that
a generic, symmetric fourth order tensor can be represented as
\begin{equation}\label{eqn:bijrs_def}
 \Bijrs = \sum_{\ga,\gb=1}^3 \Cab a_{i\alpha} a_{j\beta} a_{r\alpha}
 a_{s\beta},
\end{equation}
where $a_{ij}$ is a rotation tensor (i.e. an orthogonal matrix with
positive determinant) and $\Cab$ is a second order, symmetric tensor;
(\ref{eqn:bijrs_def}) is of course a generalization of the orthogonal
diagonalization for symmetric second order tensors.
This observation allows us to define the following algorithm:
\begin{enumerate}
\item choose a rotation tensor $a_{ij}$
\item compute with the Germano dynamic procedure the six components of
$\Cab$
\item define $\Bijrs$ using~(\ref{eqn:bijrs_def}), thereby completely
determining the subgrid flux~(\ref{eqn:tauij_aniso}).
\end{enumerate}
The anisotropic model does not prescribe how to choose the tensor
$a_{ij}$, which in principle can be any rotation tensor, possibly
varying in space and time. The values of the components $\Cab$
computed with the dynamic procedure depend on the chosen tensor,
and different choices for $a_{ij}$ result in general in different
subgrid fluxes. Many different choices have been proposed in the past,
essentially trying to identify at each position three directions
intrinsically related to the flow configuration; examples are a
vorticity aligned basis, the eigenvectors of the velocity strain rate,
or the eigenvectors of the Leonard stresses \cite{abba:2001},
\cite{abba:2003}, \cite{gibertini:2010}. In our experience, however,
the results of the simulations do not appear to have a strong dependency
on the choice of $a_{ij}$. In the present work, the components of
$a_{ij}$ are identified with those of the canonic Cartesian basis of
the three dimensional space, i.e. $a_{ij}=\delta_{ij}$, essentially
because of the simplicity of this choice and because the results
presented here are obtained for the channel flow problem, for which
the coordinate axes do identify significant directions for the
problem, namely the longitudinal, transversal and spanwise directions.

The dynamic computation of the components $\Cab$ relies on the
introduction of a test filter operator $\hat{\cdot}$. As for the
filter $\bbar{\cdot}$ introduced in section~\ref{modeq}, the precise
definition of the test filter relies on the numerical discretization
and will be given in section~\ref{dg}; here, it will suffice to point
out that the test filter is characterized by a spatial scale $\hdelta$
larger than the spatial scale $\Delta$ associated to $\bbar{\cdot}$.
The test filter is also associated to a Favre filter, denoted by
$\breve{\cdot}$, through
the Favre decomposition
\begin{equation}\label{eqn:testfavre_decomp}
 \wh{\rho \phi} = \wh{\rho} \breve{\phi},
\end{equation}
where $\phi$ stands for any of the variables in the equations
introduced in section~\ref{modeq}. Applying the test filter to the
filtered momentum equation~(\ref{filteq-momentum}) and proceeding as
in section~\ref{modeq} we arrive at
\begin{equation}
\de_t \left( \hfrho \hfu_i \right) + \de_j \left(\hfrho \hfu_i \hfu_j\right) 
+ \frac{1}{\gamma\,\MA^2}\de_i \widehat{\fp} -\frac{1}{\RE} \de_j
\widehat{\wt{\sigma}}_{ij}
= - \de_j \left( \wh{\tau}_{ij} + \Lij \right)
\label{eq:momentum-two-averages}
\end{equation}
where
\begin{equation}\label{eqn:leo_qdm}
 \Lij = \wh{\frho\fu_i\fu_j} - \hfrho\hfu_i\hfu_j
\end{equation}
is the Leonard stress tensor. Assuming now that the
model~(\ref{eqn:tauij_aniso}) can be used to represent the
right-hand-side of~(\ref{eq:momentum-two-averages}) yields
\begin{equation}
\wh{\tau}_{ij} + \Lij =
-\hfrho \hdelta^2 |\hfS| \Bijrs\breve{\fS}_{rs},
\label{eq:model-test-filter}
\end{equation}
and upon multiplying by $a_{i\alpha}a_{j\beta}$ and summing over
$i,j$, using the orthogonality of the rotation tensor,
\[
a_{i\alpha}a_{j\beta}\left( \wh{\tau}_{ij} + \Lij \right) = -\hfrho
\hdelta^2 |\hfS| \Cab a_{r\alpha}a_{s\beta} \breve{\fS}_{rs}.
\]
Substituting~(\ref{eqn:tauij_aniso}) for $\tauij$ and solving for
$\Cab$ provides the required expression
\begin{equation}\label{eq:dynamic-Cab}
 \Cab = \dfrac{ a_{i \ga} \Lij a_{j \gb} }{a_{r\ga} a_{s\gb}
  \left( \wh{\frho \Delta^2 |\fS| \fSrs} -
  \hfrho \hdelta^2 |\hfS| \breve{\fS}_{rs} \right)},
\end{equation}
and since in this work we assume $a_{ij}=\delta_{ij}$ we immediately
have
\begin{equation}
 \mathcal{C}_{ij} = \dfrac{ \Lij }{
  \left( \wh{\frho \Delta^2 |\fS| \fSij} -
  \hfrho \hdelta^2 |\hfS| \breve{\fS}_{ij} \right)}
\label{eq:Cab-expression}
\end{equation}
and
\begin{equation}\label{eqn:tauij_aniso-dij}
 \tauij =  -\frho \Delta^2 |\fS| \mathcal{C}_{ij} \fSij.
\end{equation}
The approach outlined above has some appealing features
that allow to overcome some difficulties of the Smagorinsky model.
  Firstly,  the deviatoric and isotropic parts of the subgrid stress tensor
are modelled together, without splitting of the two contributions. 
Furthermore, thanks to the anisotropy of the subgrid model,
the use of a damping function is not necessary any more to obtain correct results in the wall region,
 as it will be clearly shown by the numerical results in section \ref{results}.
We also point out that, as it will be discussed in more details in
section~\ref{dg},   the coefficients
$\mathcal{C}_{ij}$ are assumed to be averaged on each element, while they
are  not averaged  in time. This provides a local definition for such
coefficients that does not rely on the existence of any homogeneity
direction in space or quasi-stationary hypothesis in time  \cite{sagaut:2006}. Finally,
while the Smagorinsky model~(\ref{eqn:nu_smag}) is dissipative by
construction, the dynamic procedure~(\ref{eq:dynamic-Cab}) allows for
 backscattering, i.e. a positive work done by the subgrid stresses
on the mean flow. This is indeed a desirable property of the
turbulent model, yet one must ensure that the total dissipation,
resulting from both the viscous and the subgrid stresses, is positive.
This amounts to requiring
\[
\frac{1}{\RE}\fsigmaij\fSij-\tauij\fSij \geq 0,
\]
which can be ensured
by introducing a limiting coefficient in~(\ref{eqn:tauij_aniso}), so as to obtain
\begin{equation}\label{eq:limiting-beta}
\beta = 
\left\{
\begin{array}{lll}
1, & \tauij\fSij \leq 0 \\
\min\left( 1 , \frac{1}{\RE}\frac{\fsigmaij\fSij}{\tau_{kl}\fS_{kl}}
\right), &
\tauij\fSij > 0.
\end{array}
\right.
\end{equation}
Having defined the subgrid stresses, let us consider now the subgrid
terms in the energy equation, namely $\mathbf{Q}^{{\rm sgs}}$ and
$\mathbf{J}^{{\rm sgs}}$; here, we propose to treat both of them
within the same dynamic, anisotropic framework used for the subgrid
stresses. Concerning the subgrid heat flux, we let
\begin{equation}\label{eqn:Qj_aniso}
Q_i^{{\rm sgs}} = -\frho \Delta^2 |\fS| \Bir^{Q} \de_r
\fT,
\end{equation}
where $\Bir^{Q}$ is a symmetric tensor. Assuming that $\Bir^{Q}$ is
diagonal in the reference frame defined by the rotation tensor $a$ we have
\begin{equation}\label{eqn:bjstheta_def}
 \Bir^{Q} = \sum_{\ga=1}^3 \Ca^{Q} a_{i\alpha} a_{r\alpha},
\end{equation}
where the three coefficients $\Ca^Q$ can be computed locally by the
dynamic procedure. To this aim, define the temperature Leonard flux
\begin{equation}
\Li^Q = \wh{\frho \fu_i \fT} - \hfrho\hfu_i\hfT,
\label{eq:leonard-temp}
\end{equation}
apply the test filter to the filtered energy
equation~(\ref{filteq-energy}) and observe that, thanks to the
similarity hypothesis, model (\ref{eqn:Qj_aniso}) should be also
applied in the resulting equation, so that
\begin{equation}
\wh{Q}^{{\rm sgs}}_i + \Li^Q = 
 - \wh{\frho} \hdelta^2 |\hfS| \Bir^Q \de_r \hfT.
\end{equation}
Substituting~(\ref{eqn:Qj_aniso}) for $\wh{Q}^{{\rm sgs}}_i$,
multiplying by $a_{i\alpha}$, summing over $i$ and solving for
$\Ca$ yields
\begin{equation} \label{eq:dynamic-CaQ}
\Ca^Q = \frac{a_{i\alpha}\Li^Q}{a_{r\alpha}
\left( \wh{ \frho \Delta^2 |\fS| \de_r\fT } 
- \hfrho \hdelta^2 |\hfS| \de_r\hfT
\right)
}.
\end{equation}
Concerning the turbulent diffusion flux, contrary to what is done in
the Smagorinsky model, we do not neglect the term $\tau(u_i,u_k,u_k)$
in~(\ref{eqn:Jj_sgs2}), but instead adopt a scale similarity model as
in~\cite{piomelli:2000} where such term is approximated as a subgrid
kinetic energy flux 
\begin{equation}\label{eqn:kin_en_anis}
\tau(u_i,u_k,u_k) \approx
\frho \widetilde{u_i u_k u_k} - \bbar{\rho u_ku_k}\fu_i
= \frho \widetilde{u_i u_k u_k} - \frho \fu_i \widetilde{u_k u_k}.
\end{equation}
Coherently with the other subgrid terms, $\tau(u_i,u_k,u_k)$ can now
be modeled as a function of the gradient of the resolved kinetic
energy, letting
\begin{equation}\label{eqn:Jj_model}
 \tau(u_i,u_k,u_k) = - \frho \Delta^2 |\fS|
 \Bir^J \de_r\left(\dfrac{1}{2}\fu_k\fu_k\right),
\end{equation}
where
\begin{equation}\label{eqn:bjsJ_def}
 \Bir^J = \sum_{\ga=1}^3 \Ca^J a_{i\alpha} a_{r\alpha}.
\end{equation}
Introducing the kinetic energy Leonard flux
\begin{equation}
\Li^J = \hspace{-3pt}\wh{\phantom{\hspace{4pt}}\frho\fu_i\fu_k\fu_k} -
\hfrho\hfu_i\hfu_k\hfu_k
\label{eq:leonard-ekin}
\end{equation}
and proceeding exactly as for the previous terms we arrive at
\begin{equation} \label{eq:dynamic-CaJ}
\Ca^J = \frac{a_{i\alpha}\Li^J}{a_{r\alpha}
\left( \wh{ \frho \Delta^2 |\fS| \de_r
\left( \frac{1}{2}\fu_k\fu_k \right) } 
- \hfrho \hdelta^2 |\hfS| \de_r
\left( \frac{1}{2}\hfu_k\hfu_k \right)
\right)
}.
\end{equation}

\section{Discretization and filtering}
\label{dg} \indent

The equations introduced in section \ref{modeq}, including the subgrid
scale models defined in section \ref{subgrid}, will be discretized in
space by a discontinous finite element method. The DG approach
employed for the spatial discretization is analogous to that described
in \cite{giraldo:2008} and relies on the so called Local Discontinuous
Galerkin (LDG) method, see e.g. \cite{bassi:1997}, \cite{arnold:2002},
\cite{castillo:2000}, \cite{cockburn:1998b}, for the approximation of
the second order viscous terms. We provide here a concise description
of the method; a more detailed description can be found
in~\cite{maggioni:2012}.
 
In the LDG method, (\ref{filteq}) is rewritten
introducing an auxiliary variable $\bmG$, so that
\begin{equation}\label{eq:csv_auxvar}
\begin{array}{l}
\displaystyle
\de_t \bU + \Div \bF^{{\rm c}}(\bU) - \Div \bF^{{\rm v}}(\bU,\bmG) 
 + \Div \bF^{{\rm sgs}}(\bU,\bmG) = \bS \\
 \displaystyle
 \bmG - \nabla{\boldsymbol \varphi} = 0,
\end{array}
\end{equation}
where $\bU=[\frho\,,\frho\wt{\bu}^T,\frho\fe ]^T$ are the prognostic
variables, ${\boldsymbol \varphi} = [\wt{\bu}^T,\fT]^T$ are the variables
whose gradients appear in the viscous fluxes~(\ref{eq:constitutive-Favre})
as well as the turbulent ones, and $\bS$ represents the source terms.
In~(\ref{eq:csv_auxvar}), the following compact notation for the
fluxes has been used
\[
\bF^{{\rm c}} = \left[
\begin{array}{c}
 \frho\wt{\bu} \\
 \frho\wt{\bu}\otimes\wt{\bu} +
 \frac{1}{\gamma\MA^2}\bbar{p}\mathcal{I} \\
 \frho\wt{h}\wt{\bu}
\end{array}
\right],\quad
\bF^{{\rm v}} = \left[
\begin{array}{c}
 0 \\
 \frac{1}{\RE}\wt{\sigma} \\
 \frac{\gamma\MA^2}{\RE} \wt{\bu}^T \wt{\sigma}
 - \frac{1}{\kappa\RE\PR} \wt{\bq}
\end{array}
\right],
\]
 
\[
\bF^{{\rm sgs}} = \left[
\begin{array}{c}
 0 \\
 \tau \\
 \frac{1}{\kappa} \mathbf{Q}^{{\rm sgs}}
 +\frac{\gamma\MA^2}{2}\left( 
 \mathbf{J}^{{\rm sgs}} - \tau_{kk} \wt{\bu}
 \right)
\end{array}
\right],\quad
\bS = \left[
\begin{array}{c}
 0 \\
 \frho \mathbf{f} \\
\gamma\MA^2\frho \mathbf{f} \cdot \wt{\bu}
\end{array}
\right],\]
where $\tau$, $\mathbf{Q}^{{\rm sgs}}$ and $\mathbf{J}^{{\rm sgs}}$
are given by~(\ref{eqn:nu_smag}), (\ref{eqn:Qj_smag})
and~(\ref{eqn:Jj_smag}), respectively, for the Smagorinsky model,
and by~(\ref{eqn:tauij_aniso}) (including the limiting
coefficient~(\ref{eq:limiting-beta}) ), (\ref{eqn:Qj_aniso})
and~(\ref{eqn:Jj_sgs2}) together with~(\ref{eqn:kin_en_anis}) for the
anisotropic model.

The discretization is then obtained using the classical method of
lines by first introducing a space discretization and then using a
time integrator to advance in time the numerical solution. For the
time integration, we consider here the fourth order, five stage,
Strongly Stability Preserving Runge--Kutta method (SSPRK) proposed
in~\cite{spiteri:2002}. To define
the space discretization, let us first introduce a tessellation
$\mT_h$ of $\Omega$ into tetrahedral elements $K$ such that $\Omega =
\bigcup_{K\in\mT_h} K$ and $K\cap K'=\emptyset$ and define the finite
element space
\begin{equation}\label{eqn:mV_def}
\mV_h = \left\{ v_h \in L^2(\Omega): v_h|_K \in \mathbb{P}^q(K), \,
\forall K\in\mT_h \right\},
\end{equation}
where $q$ is a nonnegative integer and $\mathbb{P}^q(K)$ denotes the
space of polynomial functions of total degree at most $q$ on $K$. For
each element, the outward unit normal on $\partial K$ will be denoted
by $\bn_{\partial K}$. The numerical solution is now defined as
$(\bU_h,\bmG_h)\in(\,(\mV_h)^5\,,\,(\mV_h)^{4\times3}\,)$ such that,
$\forall K\in\mT_h$, $\forall v_h\in\mV_h$, $\forall
\br_h\in(\mV_h)^3$,
\begin{subequations}
\label{eq:DG-space-discretized}
\begin{align}
\displaystyle
\frac{d}{dt}\int_K \bU_h v_h\,d\bx
& \displaystyle
- \int_K \bF(\bU_h,\bmG_h)\cdot\nabla v_h\, d\bx
\nonumber \\[3mm]
& \displaystyle
+ \int_{\partial K} \wideparen{\bF}(\bU_h,\bmG_h)\cdot \bn_{\partial K} v_h\,
d\sigma
= \int_K \bS v_h \,d\bx,
\\[3mm] \displaystyle
\int_K \bmG_h \cdot \br_h \,d\bx
& \displaystyle
+ \int_K {\boldsymbol \varphi_h}\nabla\cdot\br_h\, d\bx
\nonumber \\[3mm]
& \displaystyle
- \int_{\partial K} \wideparen{{\boldsymbol \varphi}} \bn_{\partial
K}\cdot\br_h \, d\sigma = 0,
\end{align}
\end{subequations}
where $\bU_h=\left[ \rho_h\,,\rho_h\bu_h\,,\rho_he_h \right]^T$,
${\boldsymbol \varphi}_h=\left[ \bu_h\,,T_h \right]^T$, $\bF =
\bF^{{\rm c}}-\bF^{{\rm v}}+\bF^{{\rm sgs}}$, and $\wideparen{\bF}$,
$\wideparen{{\boldsymbol \varphi}}$ denote the so-called
 numerical fluxes. To understand the role of the numerical
fluxes, notice that~(\ref{eq:DG-space-discretized}) can be regarded as
a weak formulation of~(\ref{eqn:mV_def}) on the single element $K$
with weakly imposed boundary conditions $\wideparen{\bF}$,
$\wideparen{{\boldsymbol \varphi}}$ on $\partial K$. Hence, the
numerical fluxes are responsible for the coupling among the different
elements in $\mT_h$. There are various possible definitions of these
fluxes, and in this work we employ the Rusanov flux for
$\wideparen{\bF}$ and the centered flux for $\wideparen{{\boldsymbol
\varphi}}$; the detailed definitions can be found, for instance,
in~\cite{giraldo:2008}. To complete the definition of the space
discretization, we mention that, on each element, the unknowns are
expressed in terms of an orthogonal polynomial basis, yielding what is
commonly called a modal DG formulation, and that all the
integrals are evaluated using quadrature formulae
from~\cite{cools:2003} which are exact for polynomial orders up to
$2q$. This results in a diagonal mass matrix in the time derivative
term of~(\ref{eq:DG-space-discretized}) and simplifies the
computation of $L^2$ projections to be introduced shortly in
connection with the LES filters.

Having defined the general structure of discretized problem, we turn
now to the definition of the filter operators $\bbar{\cdot}$ and
$\wh{\cdot}$, introduced in sections~\ref{modeq}
and~\ref{subgrid:anisotropic}, respectively, with the associated Favre
decompositions. We proceed here along the lines  proposed e.g. in
\cite{collis:2002b}, \cite{collis:2002a}, \cite{vanderbos:2007},
defining the filter operators in terms of some $L^2$ projectors.
Given a subspace $\mV\subset L^2(\Omega)$, let
$\Pi_{\mV}:L^2(\Omega)\to\mV$ be the associated projector defined by
\[
\int_\Omega \Pi_{\mV}u\,v\, d\bx =
\int_\Omega u\,v\, d\bx, \qquad \forall u,v \in\mV,
\]
where the integrals are evaluated with the same quadrature rule used
in~(\ref{eq:DG-space-discretized}). For $v\in L^2(\Omega)$, the filter
$\bbar{\cdot}$ is now defined by
\begin{equation}
\bbar{v} = \Pi_{\mV_h}v,
\label{eq:filter-bar}
\end{equation}
or equivalently $\bbar{v}\in\mV_h$ such that
\begin{equation}
\int_K \bbar{v} v_h\,d\bx = \int_K v v_h\,d\bx \qquad \forall
K\in\mT_h,\quad \forall v_h\in\mV_h.
\label{eq:filter-bar-2}
\end{equation}
Notice that the application of this filter is built in the
discretization process and equivalent to it. Therefore,
once the discretization of equations (\ref{eq:csv_auxvar}) has
been performed, only $\bbar{\cdot}$ filtered
quantities are computed by the model.
To define the test filter, we  then introduce
\begin{equation}\label{eqn:mVhat_def}
\wh{\mV}_h = \left\{ v_h \in L^2(\Omega): v_h|_K \in
\mathbb{P}^{\wh{q}}(K), \, \forall K\in\mT_h \right\},
\end{equation}
where $0\leq\wh{q}<q$, and we let, for $v\in L^2(\Omega)$,
\begin{equation}
\wh{v} = \Pi_{\wh{\mV}_h}v.
\label{eq:filter-hat}
\end{equation}
By our previous identification
of the $\bbar{\cdot}$  filter and the discretization,
the quantities $\frho$, $\frho\wt{\bu}$ and $\frho\fe$
can be identified with
  $\rho_h$, $\rho_h \bu_h $ and $\rho_he_h,$ respectively.
  Therefore, they belong to $\mV_h,$ for which an orthogonal
  basis is employed by the numerical method.
  As a result, the computation of $\wh{\rho_h}$, $\wh{\rho_h\bu_h}$
and $\wh{\rho_he_h}$ is straightforward and reduces to zeroing the
last coefficients in the local expansion.
Assuming that the analytic solution is defined in some infinite
dimensional subspace of $L^2$, heuristically, $\mV_h\subset L^2$
is associated to the scales which are represented by the model, while
$\wh{\mV}_h\subset\mV_h\subset L^2$  is associated to
the spatial scales  well resolved by the numerical approximation. A
similar concept of  believable scales was introduced
in~\cite{lander:1997} in the framework of a global spectral transform
model for numerical weather prediction.

The  Favre filters associated to  (\ref{eq:filter-bar}) and (\ref{eq:filter-hat})
are defined by imposing pointwise  the conditions
(\ref{eqn:favre_decomp}--\ref{eq:constitutive-Favre})
and (\ref{eqn:testfavre_decomp}), respectively. Notice that,
as a result, for a generic quantity $\varphi$ the filtered counterpart $\wt{\varphi}$
is not, in general, a polynomial function. More specifically,
 the Favre filtered quantities are computed taking ratios of two polynomials. All
the remaining quantities in~(\ref{eqn:leo_qdm}), (\ref{eq:dynamic-Cab}),
(\ref{eq:leonard-temp}), (\ref{eq:dynamic-CaQ}),
(\ref{eq:leonard-ekin}) and~(\ref{eq:dynamic-CaJ}) where the test
filter appears are computed using~(\ref{eq:filter-hat}) and the same
quadrature rule used in~(\ref{eq:DG-space-discretized}).
We also remark that these filters
 do not commute with the differentiation operators. As
 previously remarked in section~\ref{modeq},
we neglect this error, according
to a not uncommon practice in LES
modeling~\cite{sagaut:2006}. We plan to address
this issue in more detail in a future work. An analysis of the terms resulting
from non zero commutators between differential operators and projection filters
is presented in \cite{vanderbos:2007}.

Finally, we remark that using~(\ref{eq:dynamic-Cab}), (\ref{eq:dynamic-CaQ})
and~(\ref{eq:dynamic-CaJ}), the dynamic coefficients $\Cab$,
$\Ca^Q$ and $\Ca^J$ can be computed as functions of space.
Substituting these functions directly into the subgrid dynamical models, however,
would result in diffusive terms with (possibly) highly irregular
diffusion coefficients, which would represent a serious obstacle for a
high-order numerical discretization. For this reason, the dynamic
coefficients $\Cab$, $\Ca^Q$ and $\Ca^J$ are first averaged on each
element and then used in the corresponding subgrid models. This is
similar to what is often done in the context of dynamic LES models,
where the dynamic coefficients are averaged on some homogeneity
direction, or local in space and in time~\cite{germano:1991, yang:1993, zang:1993},
with the advantage that in the present case the average is built on
the computational grid and does not require choosing any special
averaging direction. In our implementation, the dynamic coefficients
are updated at each Runge--Kutta stage; an alternative approach where
they are updated only once for each time-step or each a fixed number of time-steps could be
considered to reduce the computational cost.

Another important point is choosing the space scales $\Delta$ and
$\hdelta$ associated with the two filters~(\ref{eq:filter-bar})
and~(\ref{eq:filter-hat}). This can be done by dividing the element
diameter by the cubic root (or, in two dimensions, the square root) of
the number of degrees of freedom of $\mathbb{P}^q(K)$, for $\Delta$,
and $\mathbb{P}^{\wh{q}}(K)$, for $\hdelta$; as anticipated, this
leads to space scales which are piecewise constant on $\mT_h$. A more
precise definition, introducing a scaling coefficient which accounts
for the mesh anisotropy, is given in section~\ref{results}.

\section{Numerical results}
\label{results} \indent

In order to compare the performance of the described Smagorinsky   and  anisotropic dynamic models,
we have computed a typical  LES benchmark for compressible periodic
channel flow   at Mach numbers $\MA=0.2, 0.7, 1.5, $ respectively. 
The results obtained  are compared here with the data from the incompressible numerical simulation of 
Moser et al. (MKM) \cite{moser:1999} for $\MA=0.2$, with the simulation of Wei and Pollard (WP) \cite{wei:2011}
for $\MA=0.7$, and finally with the results presented by Coleman et al.~(CKM) \cite{coleman:1995} for the supersonic
case at $\MA=1.5$.

All the computations were performed using the {\tt FEMilaro} finite
element library~\cite{femilaro}, a FORTRAN/MPI library which,
exploiting modern FORTRAN features, aims at providing a flexible
environment for the development and testing of new finite element
formulations, and which is publicly available under GPL license.

The computational domain $\Omega^\dm$ is a box of dimensions $L^\dm_x$,
$L^\dm_y$, $L_z^\dm$ in dimensional units
 that is aligned with a 
 reference frame such that $x^\dm$ represents the streamwise axis,
 $y^\dm$ the wall normal and $z^\dm$ the spanwise axis. We also
 introduce $d^\dm = L^d_y/2$, the half height of the channel.
The reference quantities are chosen as follows
\begin{equation}
\rho_r = \rho_{{\rm b}}^\dm, \quad L_r = d^\dm, \quad
V_r = U_{{\rm b}}^\dm, \quad T_r = T^\dm_{{\rm w}},
\label{eq:channel-reference}
\end{equation}
where $\rho_{{\rm b}}^\dm$ and $U_{{\rm b}}^\dm$ are the bulk
density and the target bulk velocity, respectively, and $T^\dm_{{\rm w}}$ is the wall
temperature. In dimentionless units we let $L_x=4\pi$, $L_y=2$ and
$L_z=2\pi$ for all the computations, except the cases with $\MA=0.2$
where we choose $L_x=2\pi$; the resulting domain is thus
$\Omega=[0\,,4\pi]\times[-1\,,1]\times[0\,,2\pi]$,
or $\Omega=[0\,,2\pi]\times[-1\,,1]\times[0\,,2\pi]$ for $\MA=0.2$.
Isothermal, no-slip boundary conditions are imposed for $y=\pm1$, i.e.
$T=1$ and $\bu = 0$, while periodic conditions are applied in the
streamwise and spanwise directions. The initial condition is
represented by a laminar Poiseuille profile $u_x =
\frac{3}{4}(1-y^2)$, with $\rho=1$ and $T=1$. A random perturbation of
amplitude $a=0.1$
is added to the initial velocity, while no perturbations are added to
$\rho$ and $T$. The perturbation of the $(i+1)$-th velocity component
is evaluated at each quadrature node by scaling the $i$-th coordinate
of the node to obtain $\xi^{(0)} \in (0\,,1)$, computing 20 iterations
of the logistic map $\xi^{(k+1)} = 3.999\,\xi^{(k)}(1-\xi^{(k)})$ and
projecting the resulting values, which turn out to be uncorrelated in
space, on the local polynomial space; this provides a simple,
deterministic and portable way to define a random perturbation of the
velocity with zero divergence. The value $U^\dm_{{\rm b}}$ is, by
definition, the desired
bulk velocity; the flow velocity, however, is the result of the
balance between the external forcing and the dissipations, and can not
be easily fixed \emph{a priori}. To ensure that the obtained bulk
velocity coincides with the prescribed value, as well as
 to preserve the
homogeneity of the flow in the directions parallel to the wall,
a body force uniform in space is included along the
streamwise direction, defined by
\begin{equation}\label{forcing}
f_x(t) = - \frac{1}{\rho_{{\rm b}}} \left[ \alpha_1 \left( Q(t) - Q_0 \right) + \alpha_2 \int_0^t
\left( Q(s)- Q_0 \right) ds \right],
\end{equation}
where $Q(t)=\int_\Omega \rho(t,\bx) u_x(t,\bx)d\bx /L_x$ is the instantaneous flow
rate and $Q_0=L_yL_z$ is the flow rate corresponding to the desired bulk
velocity.
A sufficiently rapid convergence toward the value $Q_0$
has been observed by taking $\alpha_1=0.1$, $\alpha_2=0.5$.
The bulk Reynolds and Mach numbers are defined as
\begin{equation}
\RE_{{\rm b}} = \frac{\rho_{{\rm b}}^\dm U_{{\rm b}}^\dm
d^\dm}{\mu_{{\rm w}}^\dm}, \qquad
\MA_{{\rm b}} = \frac{U_{{\rm b}}^\dm}{\sqrt{\gamma R T_{{\rm
w}}^\dm}},
\label{eq:channel-ReMa}
\end{equation}
where $\mu_{{\rm w}}^\dm$ is the viscosity at the wall.

The computational mesh employed is obtained by a structured 
mesh with $N_x=16$ ($N_x=8$ for $\MA=0.2$), $N_y=16$, $N_z=12$ hexahedra in the $x, y, z$ directions, respectively, 
each of which is then split into $N_t=6$ tetrahedral finite elements.
 While uniform in the $x,z$ directions, the hexahedral mesh  is not
 uniform in the $y$ direction, where
the $y=const$ planes are given by
\begin{equation}\label{eqn:stretch}
 y_j = -\dfrac{\tanh\left(\w\left(1-2j/N_y\right)\right)}
              {\tanh\left(\w\right)} \qquad \mbox{for} \qquad j=0, \ldots, N_y .
\end{equation}
The value of the parameter $\w $ is chosen by fixing the position
$y_1$ for the face closest to the wall so that the laminar sublayer is
well resolved. For each tetrahedral element $K,$ we will then denote
by $\Delta^{(l)}_K, \ l=1,2,3 $ the dimensions of the hexahedron from
which the element was obtained in the $x,$ $y,$ and $z $ coordinate
directions, respectively. The polynomial degrees for $\mV_h$ and
$\wh{\mV}_h$ are $q=4$ and $\wh{q}=2$, respectively.
For the basis functions of degree 4, 
 $N_q=35$ number of degrees of freedom were employed in each element, while for  the basis functions of degree 2
 $N_{\wh{q}}=10$ degrees of freedom were employed. 
As a result, the grid spacing is given in the homogeneity directions by
$$\Delta_x=\frac{L_x}{N_x \sqrt[3]{N_t N_q}} \quad \quad  
\Delta_z=\frac{L_z }{N_z \sqrt[3]{N_t N_q}}.$$
The grid filter scale $\Delta(K) $ can then be estimated 
as suggested by \cite{lilly:1993} for strongly anisotropic grids. 
For each element $K,$ we define
$$
\Delta_{max}(K)  = \max_i \Delta^{(i)}(K)  \qquad  
 a_l = \frac{\Delta^{(l)}(K)}{\Delta_{max} (K)}  \qquad   a_k = \frac{\Delta^{(k)}(K)}{\Delta_{max} (K)}  
$$
where $l$ and $k$ are the directions in which the maximum is not attained, and
\begin{subequations}
\begin{align}
 f &= \cosh \sqrt{\dfrac{4}{27}\left[ (\ln a_l)^2 - \ln a_l \ln a_k + (\ln a_k)^2 \right]} \\
 \Delta (K)  &=  \left( \frac{\prod_{i=1}^3\Delta^{(i)}(K)}{N_q} \right)^{1/3} f. \label{eqn:f-lilly}
\end{align}
\end{subequations}
The test filter scale  $\hdelta (K) $ is defined analogously, only
replacing  $N_q$ by $N_{\wh{q}} $ in the previous definitions.
The parameters for the three cases considered here  and for the comparison test cases presented in literature are 
summarized in Table \ref{tab:parameter}.
\begin{table}
{\footnotesize 
\begin{tabular}[hbtp]{|l|c|c|c|c|c|c|}
\hline
& Moser  & Wei and  & Coleman & Present  & Present & Present \\
& et al. &  Pollard & et al. & Ma=0.2 &  Ma=0.7 & Ma=1.5 \\
 & (MKM) & (WP) & (CKM) & (Ma02) & (Ma07) & (Ma15) \\
 \hline
$\MA_{{\rm b}}$ & \textemdash & 0.7 & 1.5 & 0.2 & 0.7 & 1.5 \\
 \hline
$\RE_{{\rm b}}$ & 2800 & 2795 & 3000 & 2800 & 2795 & 3000  \\
 \hline
 $L_x$ & $4\pi$ &  12 & $4\pi$ & $2\pi$ & $4\pi$ & $4\pi$ \\
 \hline
 $L_z$ & $\frac{4}{3}\pi$ & 6 & $\frac{4}{3}\pi$ &  $\frac{4}{3}\pi$ &  $\frac{4}{3}\pi$ & $\frac{4}{3}\pi$ \\
 \hline
 $\Delta_x^+$ & 17.7 &  4.89 & 19 & 23 & 24  & 29 \\
 \hline
 $\Delta_z+$ & 5.9 & 4.89 & 12 & 10 &  11 & 13 \\
 \hline
 $\Delta^+_{y_{min}}/\Delta^+_{y_{max}}$ &  0.05/4.4 & 0.19/2.89 & 0.1/5.9 & 0.65/7.9 & 0.67/8.2 & 0.8/9.5 \\
\hline
\end{tabular}
}
\caption{Parameters of simulations and reference test cases.}
\label{tab:parameter}
\end{table}
For the present case, the grid spacing $\Delta_x^+$, $\Delta_y^+$, $\Delta_z^+$ in wall unit have been estimated
{\em a posteriori} as
\begin{equation}
 \Delta_x^+ = \Delta_x \RE_\tau, \quad \Delta_y^+ = \Delta_y \RE_\tau,
 \quad \Delta_z^+=\Delta_z \RE_\tau,
 \nonumber
\end{equation}
where $\RE_\tau$ is the skin friction Reynolds number otained by the simulations 
and reported in Table \ref{tab:mean_var}.

For the Smagorinsky-type model, a test with $C_I=0.01$ seemed to enhance the dissipative behaviour of the
model, so that all the results  presented  in the following have been computed with $C_I=0$, as in 
\cite{erlebacher:1992} and \cite{lenormand:2000} where the isotropic contribution is neglected.

After the statistical steady state was reached at time $t_{{\rm st}}$, the simulations
were continued for a dimensionless time $t_{{\rm av}}$ at least equal to $60$ non dimensional time units
to compute all the relevant statistics and verifying time invariance of mean profiles.
The statistics are now computed averaging on the element faces
parallel to the walls, introducing, for a generic quantity $\varphi$, the space-time average
\begin{equation}
\begin{array}{ll}
\displaystyle
<\varphi>(|y|) = \frac{1}{2t_{{\rm av}}L_xL_z} \\[3mm]
\displaystyle \qquad
\int_{t_{{\rm st}}}^{t_{{\rm st}}+t_{{\rm av}}} \int_0^{L_x}
\int_0^{L_z} \left(\varphi(t,x,-|y|,z) +
\varphi(t,x,|y|,z)\right)\,dz\,dx\,dt.
\end{array}
\label{eq:channel-average}
\end{equation}

In Table \ref{tab:mean_var} the mean flow quantities at the wall and at the channel centerline,
denoted by the subscripts ${\rm w}$  and ${\rm c}$, respectively, are compared 
with the reference DNS results. 
\begin{table}
{\footnotesize 
\begin{tabular}[hbtp]{|l|c|c|c|c|c|c|c|c|}
\hline
& $\tau_{{\rm w}}$ & $Re_\tau$ & $u^\dm_\tau/U^\dm_{{\rm b}}$ &
$\rho^\dm_{{\rm w}}/\rho^\dm_{{\rm b}}$ &
$U_{{\rm c}}^\dm/U^\dm_{{\rm b}}$ & $\rho^\dm_{{\rm c}}/\rho^\dm_{{\rm b}}$ &
$\rho_{{\rm c}}^\dm/\rho_{{\rm w}}^\dm$ & $T_{{\rm c}}^\dm/T^\dm_{{\rm w}}$ \\
\hline
MKM          & $11.21$ & $178$ & $0.06357$ & \textemdash & $1.1672$ & \textemdash & \textemdash & \textemdash \\
\hline  
Anis. Ma02 & $10.08$ & $169$ & $0.05995$ & $1.0037$    & $1.1355$ & $0.9998$    & $0.9962$    & $0.9973$ \\
\hline
Smag. Ma02 & $9.98$ & $167$ & $0.05964$ & $1.0037$    & $1.1613$ & $0.9998$    & $0.9962$    & $1.005$  \\
\hline
WP         & $12.38$  & $186$ & $0.06184$ & $1.1076$    & $1.1636$ & $0.9949$    & $0.9246$    & $1.0863$ \\
\hline
Anis. Ma07 & $10.22$ & $169$ & $0.057$   & $1.0649$    & $1.1613$ & $0.9961$    & $0.9353$    & $1.071$ \\
\hline
Smag. Ma07 & $9.20$ & $160$ & $0.0502$  & $1.0624$    & $1.1691$ & $0.9959$    & $0.9374$    & $1.070$ \\
\hline
CKM        & $12.12$ & $222$ & $0.0545$  & $1.3578$    & $1.164$  & $0.9817$    & $0.723$     & $1.378$ \\
\hline
Anis. Ma15 & $11.30$ & $209$ & $0.05404$ & $1.2898$    & $1.1513$ & $0.983$     & $0.7621$    & $1.335$ \\
\hline
Smag. Ma15 & $9.94$ & $194$ & $0.05122$ & $1.2632$    & $1.1744$ & $0.9845$    & $0.7794$    & $1.313$ \\
\hline
\end{tabular}
}
\caption{Mean flow quantities for all the numerical experiments.}
\label{tab:mean_var}
\end{table}
In the Ma02 simulations, the constant density and temperature conditions of the incompressible MKM DNS 
are recovered with an error in the order of $3\tcperthousand $ at most. 
The wall shear stress $\tau_{{\rm w}}=\mu_{{\rm
w}}(\frac{\partial}{\partial y}<u>)_{{\rm w}}$ is the most sensitive quantity and is always
underestimated. The wall stress relative errors range between $6\div25\%$, where the
larger values are obtained with the Smagorinsky model.
The Reynolds number $\RE_\tau=\sqrt{\rho_{{\rm w}} Re_{{\rm b}} (\frac{\partial}{\partial y}<u>)_{{\rm w}}}$ and the skin-friction velocity
$u_\tau=\RE_\tau/(\RE_{{\rm b}} \rho_{{\rm w}})$ are affected by the wall shear stress error and by the fact that the density 
$\rho_{{\rm w}}$ at the wall   is always underpredicted. On the other hand, at the center of the channel 
density values are higher than the reference ones and, coherently,  temperature values are lower. 
The mean velocity at the centerline is always underestimated,
 except for the compressible cases computed with the Smagorinsky  model.
 The overprediction of this quantity by the Smagorinsky model is probably related
 to its difficulties in connecting properly the wall region to the
  the logarithmic layer.
Looking at the mean quantities, for all  Mach number values and all indicators considered, 
the anisotropic model performs as well as or better than the Smagorinsky model, especially in the wall region.

\begin{figure}
\centering
\begin{subfigure}[]{\label{fig:rhomean_Ma07}
    \includegraphics[width=0.48\textwidth]{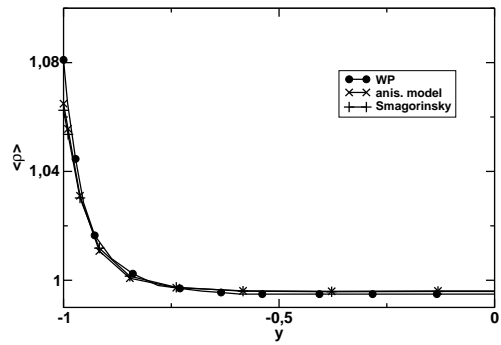}}    
\end{subfigure} 
\begin{subfigure}[]{\label{fig:rhomean_Ma15}
    \includegraphics[width=0.48\textwidth]{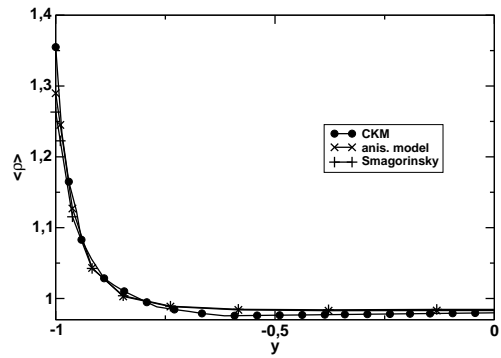}}
\end{subfigure} 
\caption{Mean density profiles at $Ma=0.7$ 
(Fig.\ref{fig:rhomean_Ma07}) and $Ma=1.5$ (Fig.\ref{fig:rhomean_Ma15}).}
\label{fig:rhomean}
\end{figure}


\begin{figure}
\centering
\begin{subfigure}[]{\label{fig:Tmean_Ma07}
    \includegraphics[width=0.48\textwidth]{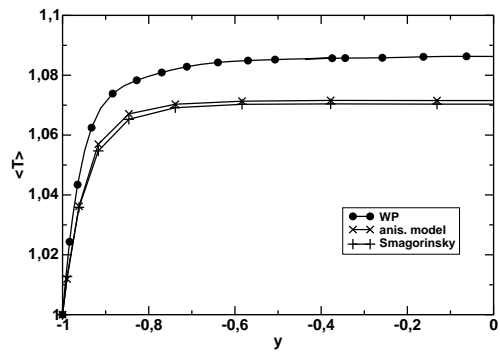}}    
\end{subfigure} 
\begin{subfigure}[]{\label{fig:Tmean_Ma15}
    \includegraphics[width=0.48\textwidth]{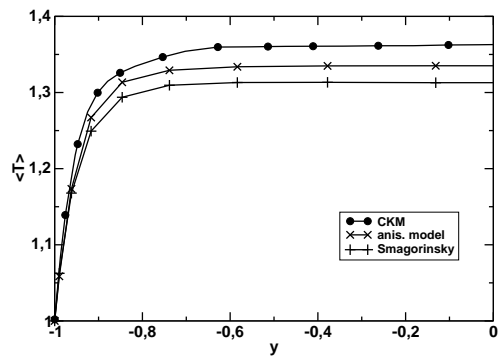}}
\end{subfigure} 
\caption{Mean temperature  profiles at $Ma=0.7$ 
(Fig.\ref{fig:Tmean_Ma07}) and $Ma=1.5$ (Fig.\ref{fig:Tmean_Ma15}).}
\label{fig:Tmean}
\end{figure}

\begin{figure}[htbc]
\centering
\begin{subfigure}[]{\label{fig:umean_Ma02}
      \includegraphics[width=0.55\textwidth]{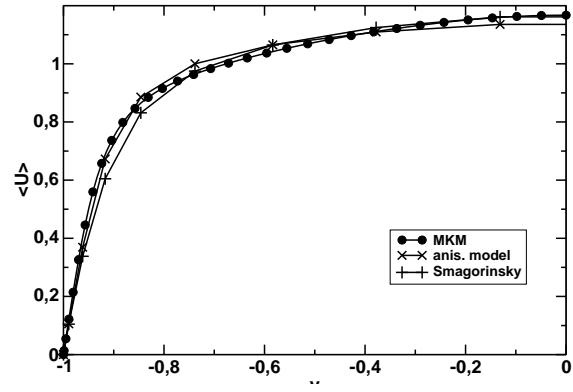}}
\end{subfigure}
\begin{subfigure}[]{\label{fig:umean_Ma07}
      \includegraphics[width=0.55\textwidth]{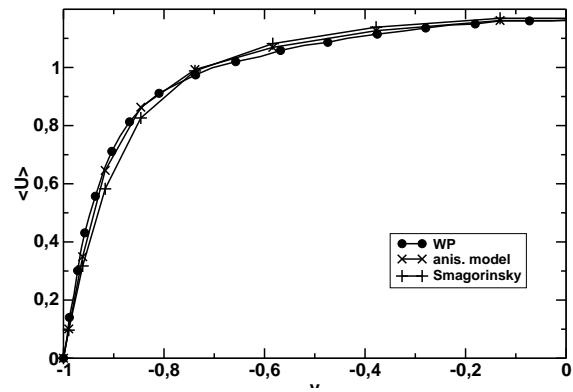}}
\end{subfigure} 
\begin{subfigure}[]{\label{fig:umean_Ma15}
      \includegraphics[width=0.55\textwidth]{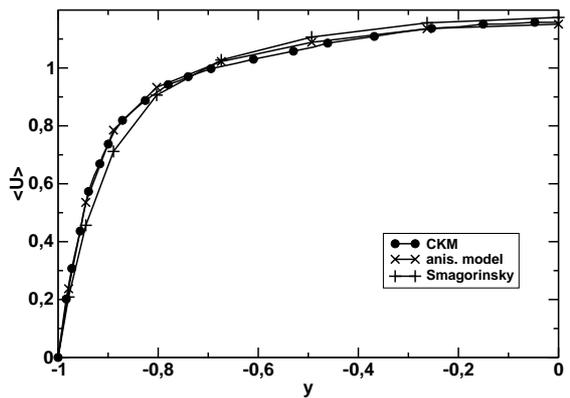}}
\end{subfigure} 
\caption{Mean streamwise velocity  profiles at  (a) $\MA=0.2,$ (b)
$\MA=0.7$ and
(c) $\MA=1.5.$}
\label{fig:umean}
\end{figure}

 Figure \ref{fig:rhomean} confirms the mean density values reported in Table \ref{tab:mean_var}. 
The excess in the density profiles at the channel center is related to the temperature values lower than the DNS 
ones far from the wall (see Figure \ref{fig:Tmean}). In spite of this, in Figure \ref{fig:Tmean} the mean temperature 
profiles demonstrate the improvement due to the modeling of subgrid
terms in the energy equations with the anisotropic model with respect to the Smagorinsky one, especially in the 
supersonic case.  Figure \ref{fig:umean} shows instead 
the mean velocity profiles.
 It is apparent that the anisotropic model approximates better
 the DNS results close to the wall.  
 
\begin{figure}[htbc]
\centering
\includegraphics[width=0.55\textwidth]{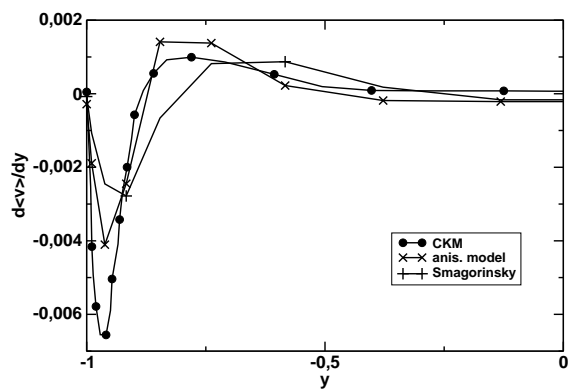}
\caption{Mean dilatation profiles at $\MA=1.5.$}
\label{fig:dvdy}
\end{figure}

Figure \ref{fig:dvdy} shows the mean profile of the non-solenoidal
term\linebreak
 $\frac{\partial}{\partial y}<\nolinebreak v\nolinebreak>$
 in the supersonic case. With the anisotropic model, the
 compression near the wall is underestimated, but the peak position is well captured, while this is not the case for the 
Smagorinsky model. At the center of the channel, while for the DNS a small dilatation is present, for the
LES a small compression is probably necessary to compensate the excess of dilatation taking place
in the buffer layer between the wall region and the logarithmic layers.

\begin{figure}[htbc]
\centering
\begin{subfigure}[]{\label{fig:rmsv_Ma02}
    \includegraphics[width=0.55\textwidth]{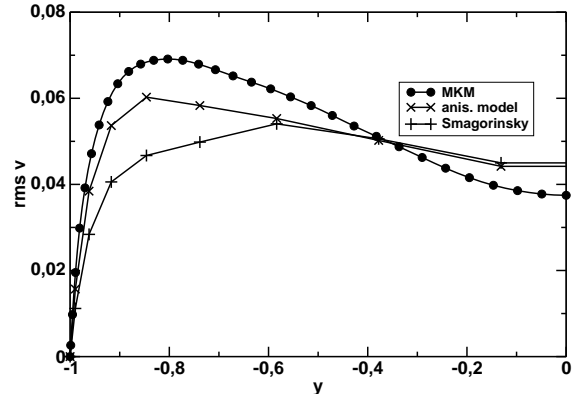}}
\end{subfigure} 
\begin{subfigure}[]{\label{fig:rmsv_Ma07}
    \includegraphics[width=0.55\textwidth]{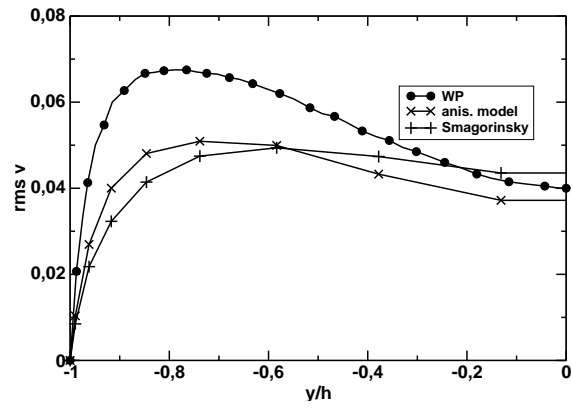}}
\end{subfigure} 
\begin{subfigure}[]{\label{fig:rmsv_Ma15}
    \includegraphics[width=0.55\textwidth]{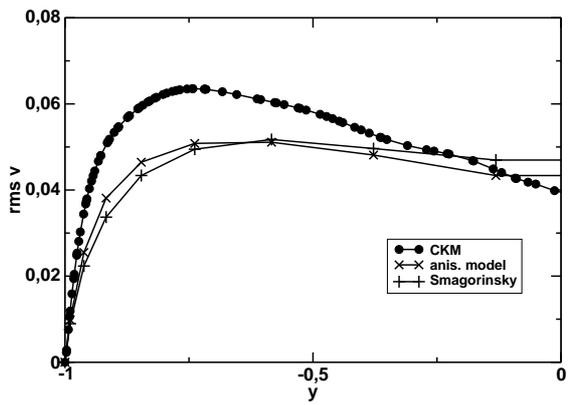}}
\end{subfigure} 
\caption{Root mean square profiles of the wall normal velocity
component at (a) $\MA=0.2$,
(b) $\MA=0.7$ and (c) $\MA=1.5$.}
\label{fig:rmsv}
\end{figure}

In figures \ref{fig:rmsu}-\ref{fig:rmsw},  the root mean square values of the resolved 
velocity fluctuations are displayed.
Figure \ref{fig:rmsu} for the streamwise turbulence intensity shows that the Smagorinsky model presents 
different behavior depending on the Mach number. In the incompressible 
limit, the dissipative character of the Smagorinsky model leads to an overprediction of the 
streamwise turbulence intensity in the wall region, as often verified
in other LES experiments (see for instance
\cite{lenormand:2000}). We recall that these quantities represent the resolved contributions only, so that their 
overestimation with respect to the DNS value is an undesired result.
In the compressible simulations, instead, the Smagorinsky model  predicts well  the value of the peak but not the
position. Moreover, the fluctuations 
are larger than those of  the DNS   far from the wall. On the other hand, the fluctuation peak 
is always well captured by the anisotropic model. Furthermore,   while
for $\MA=0.2$ the streamwise intensities are overpredicted by the anisotropic model in the center region, in the other tests they are well estimated.

 The fluctuations of the  velocity components normal to the wall (Figure \ref{fig:rmsv}) and  spanwise (Figure \ref{fig:rmsw}) in the wall region are  underestimated by both models 
 with respect to the DNS values, although we recall that these 
 are the resolved contribution only. 
 For these components the difference between the Smagorinsky and the anisotropic model become less evident
as the Mach number increases, but   the anisotropic model always performs better.
 As it usually happens in LES experiments, 
  lower values of the wall normal components obtained with the Smagorinsky model are
 associated to the overprediction of the streamwise fluctuations. In the centerline region, 
where the turbulence
presents a more isotropic character, the anisotropic model tends to 
recover the Smagorinsky model results, especially  at $\MA=0.2$.   

\begin{figure}[htbc]
\centering
\begin{subfigure}[]{
    \includegraphics[width=0.55\textwidth]{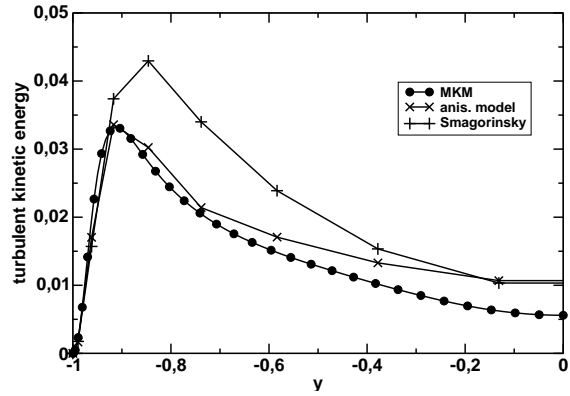}}
\end{subfigure}
\begin{subfigure}[]{ 
    \includegraphics[width=0.55\textwidth]{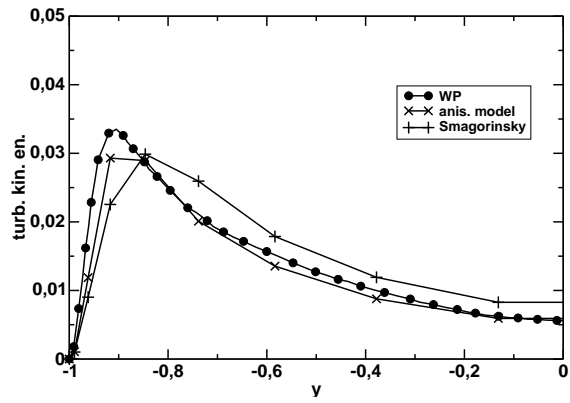}}
\end{subfigure}
\begin{subfigure}[]{ 
    \includegraphics[width=0.55\textwidth]{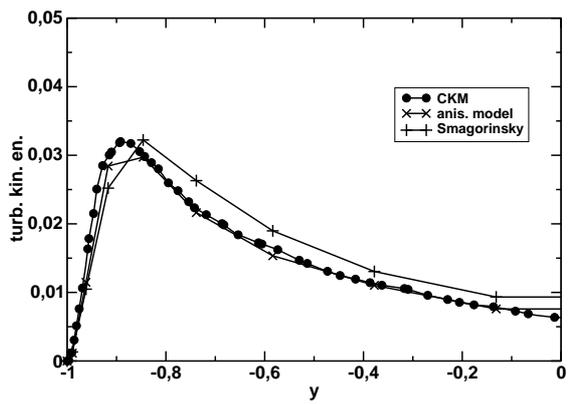}}
\end{subfigure}
\caption{Total modelled+resolved turbulent kinetic energy at (a) $\MA=0.2$,
(b) $\MA=0.7$ and (c) $\MA=1.5$.}
\label{fig:turb_kinen}
\end{figure}

In Figure \ref{fig:turb_kinen} results for the total
(modelled plus resolved) turbulent kinetic energy are displayed. For the Smagorinsky model, this corresponds 
to the resolved turbulent kinetic energy, since the isotropic part of the subgrid stresses is neglected.
It can be observed that also for this quantity the DNS results  are very well reproduced by the anisotropic model.

\begin{figure}[htbc]
\centering
\begin{subfigure}[]{
    \includegraphics[width=0.55\textwidth]{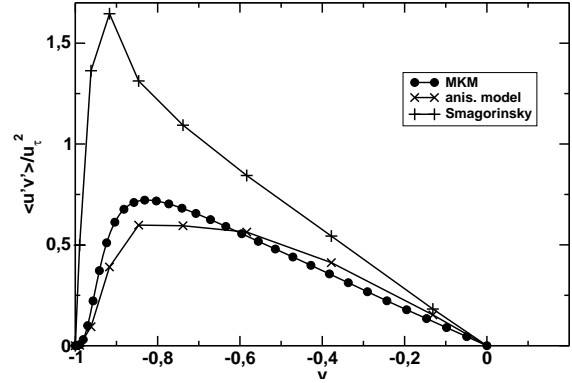}}
\end{subfigure}
\begin{subfigure}[]{ 
    \includegraphics[width=0.55\textwidth]{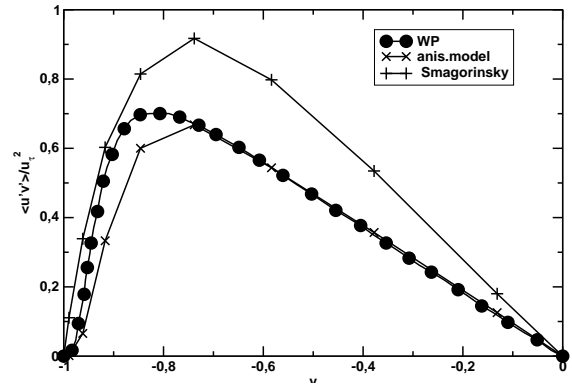}}
\end{subfigure}
\begin{subfigure}[]{ 
    \includegraphics[width=0.55\textwidth]{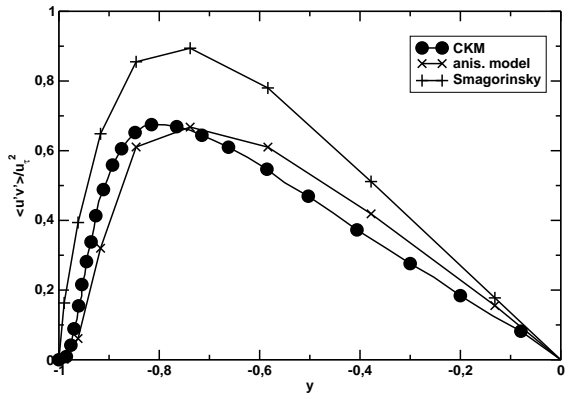}}
\end{subfigure}
\caption{Total modelled+resolved turbulent wall normal shear stress at
(a) $\MA=0.2$,
(b) $\MA=0.7$ and (c) $\MA=1.5$. The stress is normalized by the corresponding $u_\tau$ wall 
friction velocity obtained by the simulation.}
\label{fig:uv}
\end{figure}

Since during the simulations a constant mass flow is imposed, the wall shear stress $\tau_{{\rm w}}$ 
can differ from the expected DNS  value (see Table \ref{tab:mean_var}) and relevant differences affect also 
the wall normal turbulent shear stress (modeled + resolved) reported in Figure \ref{fig:uv}. 
Here, the stress is rescaled by the corresponding $u_\tau$ wall friction velocity obtained by each simulation.
In spite of the application of  the damping function, the Smagorinsky model does not 
present the correct trend at the wall and the shear 
stress is overestimated. This behaviour is probably the cause of the underprediction of the mean velocity profile
in the wall region and of  its difficulties in connecting properly the wall region to the
  the logarithmic layer. On the other hand,
 the anisotropic model 
is in quite good agreement with the DNS results for  simulations at all  Mach numbers.
 
\begin{figure}[htbc]
\centering
\begin{subfigure}[]{\label{fig:rmsu_Ma02}
    \includegraphics[width=0.55\textwidth]{figure/rmsu_Ma02.eps}}
\end{subfigure} 
\begin{subfigure}[]{\label{fig:rmsu_Ma07}
    \includegraphics[width=0.55\textwidth]{figure/rmsu_Ma07.eps}}
\end{subfigure} 
\begin{subfigure}[]{\label{fig:rmsu_Ma15}
    \includegraphics[width=0.55\textwidth]{figure/rmsu_Ma15.eps}}
\end{subfigure} 
\caption{Root mean square profiles of the streamwise velocity
component at (a) $\MA=0.2$,
(b) $\MA=0.7$ and (c)  $\MA=1.5$.}
\label{fig:rmsu}
\end{figure}

\begin{figure}[htbc]
\centering
\begin{subfigure}[]{\label{fig:rmsw_Ma02}
    \includegraphics[width=0.55\textwidth]{figure/rmsw_Ma02.eps}}
\end{subfigure} 
\begin{subfigure}[]{\label{fig:rmsw_Ma07}
    \includegraphics[width=0.55\textwidth]{figure/rmsw_Ma07.eps}}
\end{subfigure} 
\begin{subfigure}[]{\label{fig:rmsw_Ma15}
    \includegraphics[width=0.55\textwidth]{figure/rmsw_Ma15.eps}}
\end{subfigure} 
\caption{Root mean square profiles of the spanwise velocity components
at (a) $\MA=0.2$,
(b) $\MA=0.7$ and (c) $\MA=1.5$.}
\label{fig:rmsw}
\end{figure}

%
 
\section{Conclusions and future perspectives}
\label{conclu} \indent

 We have investigated the potential benefits resulting from the application  of
the anisotropic dynamic model  \cite{abba:2003} in the context of a high order DG model.
This approach contrasts with other attempts at implementing LES in a DG framework,
in which  only Smagorinsky closures have been applied so far.
Furthermore, the hierarchical nature of the DG finite element basis was exploited to
implement the LES grid and test filters via projections on the finite dimensional
subspaces that define the numerical approximation, along the lines of similar proposals
in the VMS framework. A comparison with the DNS
experiment results reported in  \cite{coleman:1995}, \cite{moser:1999} and \cite{wei:2011} has been carried out. 
The results of the comparison show a clear improvement in the prediction of several key features of the flow with respect to the Smagorinsky closure implemented in the same framework.
 The proposed approach appears to lead to significant improvements both in the
low and high  Mach number regimes. On this basis, we plan to investigate 
 further extensions of this approach to flows in presence of gravity, 
with the goal of improving the turbulence 
models for applications to environmental stratified flows. Furthermore, the numerical framework that has
been  validated by the comparison reported in this paper 
will be employed for the assessment of the proposal presented in \cite{germano:2014} for the
extension of the eddy viscosity model to compressible flows.

 \section*{Acknowledgements} 
 Part of the results have been already presented in the Master thesis in Aerospace
Engineering of A. Maggioni, prepared at Politecnico di Milano under the supervision of some of the authors.
We would like to thank A. Maggioni for the first implementation
of the dynamic models we have employed in this work. We are also very grateful to M.Germano for
several useful discussions on the topics studied in this paper.
The present research has been carried out with financial 
support by Regione Lombardia and by the Italian Ministry of Research and Education 
in the framework of the  PRIN 2008 project  
{\it Analisi e sviluppo di metodi numerici avanzati per Equazioni 
alle Derivate Parziali.} We acknowledge that the results of this research have been achieved using the 
computational resources made available    at CINECA (Italy) by the high performance computing projects
ISCRA-C HP10CAM1FM and HP10CVWE4N.

\bibliographystyle{plain}
\bibliography{dg_les}

\begin{thebibliography}{10}

\bibitem{abba:2001}
A.~Abb\`a, C.~Cercignani, G.~Picarella, and L.~Valdettaro.
\newblock A 3{D} turbulent boundary layer test for {LES} models.
\newblock In {\em Computational Fluid Dynamics 2000}, 2001.

\bibitem{abba:2003}
A.~Abb\`a, C.~Cercignani, and L.~Valdettaro.
\newblock {Analysis of {S}ubgrid {S}cale {M}odels}.
\newblock {\em Computer and Mathematics with Applications}, 46:521--535, 2003.

\bibitem{arnold:2002}
D.N. Arnold, F.~Brezzi, B.~Cockburn, and L.D. Marini.
\newblock Unified analysis of {D}iscontinuous {G}alerkin methods for elliptic
  problems.
\newblock {\em SIAM Journal of Numerical Analysis}, 39:1749--1779, 2002.

\bibitem{bassi:1997}
F.~Bassi and S.~Rebay.
\newblock {A High Order Accurate Discontinuous Finite Element Method for the
  Numerical Solution of the Compressible {N}avier-{S}tokes Equations}.
\newblock {\em Journal of Computational Physics}, 131:267--279, 1997.

\bibitem{castillo:2000}
P.~Castillo, B.~Cockburn, I.~Perugia, and D.~Sch{\"o}tzau.
\newblock An a priori analysis of the {L}ocal {D}iscontinuous {G}alerkin method
  for elliptic problems.
\newblock {\em SIAM Journal of Numerical Analysis}, 38:1676--1706, 2000.

\bibitem{cockburn:1998b}
B.~Cockburn and C.~Shu.
\newblock {The {L}ocal {D}iscontinuous {G}alerkin Method for Time-Dependent
  Convection-Diffusion Systems}.
\newblock {\em SIAM Journal of Numerical Analysis}, 35:2440--2463, 1998.

\bibitem{coleman:1995}
G.N. Coleman, J.~Kim, and R.D. Moser.
\newblock A numerical study of turbulent supersonic isothermal-wall channel
  flow.
\newblock {\em Journal of Fluid Mechanics}, 305:159--183, 1995.

\bibitem{collis:2002b}
S.~S. Collis.
\newblock Discontinuous {G}alerkin methods for turbulence simulation.
\newblock In {\em Proceedings of the 2002 Center for Turbulence Research Summer
  Program}, pages 155--167, 2002.

\bibitem{collis:2002a}
S.~S. Collis and Y.~Chang.
\newblock The {DG}/{VMS} method for unified turbulence simulation.
\newblock {\em AIAA paper}, 3124:24--27, 2002.

\bibitem{cools:2003}
R.~Cools.
\newblock {An {E}ncyclopaedia of {C}ubature {F}ormulas}.
\newblock {\em Journal of Complexity}, 19:445--453, 2003.

\bibitem{vanderbos:2007}
F.van der Bos, J.J.W. van~der Vegt, and B.J. Geurts.
\newblock A multi-scale formulation for compressible turbulent flows suitable
  for general {V}ariational discretization techniques.
\newblock {\em Computer Methods in Applied Mechanics and Engineering},
  196:2863--2875, 2007.

\bibitem{eidson:1985}
T.M. Eidson.
\newblock Numerical simulation of turbulent {R}ayleigh-{B}\'{e}nard problem
  using subgrid modeling.
\newblock {\em Journal of Fluid Mechanics}, 158:245--268, 1985.

\bibitem{erlebacher:1992}
G.~Erlebacher, M.Y. Hussaini, C.G. Speziale, and T.A. Zang.
\newblock Large {E}ddy {S}imulation of compressible turbulent flows.
\newblock {\em Journal of Fluid Mechanics}, 238:155--185, 1992.

\bibitem{femilaro}
{FEM}ilaro, a finite element toolbox.
\newblock \url{https://code.google.com/p/femilaro/}.
\newblock Available under GNU GPL v3.

\bibitem{garcia:2014}
F.~Garcia, L.~Bonaventura, M.~Net, and J.~S\'anchez.
\newblock Exponential versus {IMEX} high-order time integrators for thermal
  convection in rotating spherical shells.
\newblock {\em Journal of Computational Physics}, 264:41--54, 2014.

\bibitem{garnier:2009}
E.~Garnier, N.~Adams, and P.~Sagaut.
\newblock {\em Large {E}ddy {S}imulation for {C}ompressible {F}lows}.
\newblock Springer Verlag, 2009.

\bibitem{germano:1992}
M.~Germano.
\newblock Turbulence: the filtering approach.
\newblock {\em Journal of Fluid Mechanics}, 238:325--336, 1992.

\bibitem{germano:1991}
M.~Germano, U.~Piomelli, P.~Moin, and W.~H. Cabot.
\newblock {A {D}ynamic {S}ubgrid-{S}cale {E}ddy {V}iscosity {M}odel}.
\newblock {\em Physics of Fluids}, 3(7):1760--1765, 1991.

\bibitem{gibertini:2010}
G.~Gibertini, A.~Abb\`a, F.~Auteri, and M.~Belan.
\newblock Flow around two in-tandem flat plates: Measurements and computations
  comparison.
\newblock In {\em 5th International Conference on Vortex Flows and Vortex
  Models (ICVFM2010)}, 2010.

\bibitem{giraldo:2008}
F.X. Giraldo and M.~Restelli.
\newblock A study of spectral element and discontinuous {G}alerkin methods for
  the {N}avier-{S}tokes equations in nonhydrostatic mesoscale atmospheric
  modeling: equation sets and test cases.
\newblock {\em Journal of Computational Physics}, 227:3849--3877, 2008.

\bibitem{giraldo:2010}
F.X. Giraldo, M.~Restelli, and M.~L\"auter.
\newblock Semi-implicit formulations of the {N}avier-{S}tokes equations:
  application to nonhydrostatic atmospheric modeling.
\newblock {\em SIAM Journal of Scientific Computing}, 32:3394--3425, 2010.

\bibitem{hughes:1998}
T.J.R. Hughes, G.R. Feijoo, L.~Mazzei, and J.B. Quincy.
\newblock The {V}ariational {M}ultiscale method-a paradigm for computational
  mechanics.
\newblock {\em Computer Methods in Applied Mechanics and Engineering},
  166:3--24, 1998.

\bibitem{hughes:2000}
T.J.R. Hughes, L.~Mazzei, and K.~Jansen.
\newblock Large {E}ddy {S}imulation and the {V}ariational {M}ultiscale method.
\newblock {\em Computing and Visualization in Science}, 3:47--59, 2000.

\bibitem{hughes:2001a}
T.J.R. Hughes, L.~Mazzei, A.A. Oberai, and A.A. Wray.
\newblock The {M}ultiscale formulation of {L}arge {E}ddy {S}imulation: Decay of
  homogeneous isotropic turbulence.
\newblock {\em Physics of Fluids}, 13:505--512, 2001.

\bibitem{hughes:2001b}
T.J.R. Hughes, A.A. Oberai, and L.~Mazzei.
\newblock Large{E}ddy {S}imulation of turbulent channel flows by the
  {V}ariational {M}ultiscale method.
\newblock {\em Physics of Fluids}, 13:1784--1799, 2001.

\bibitem{hughes:2004}
T.J.R. Hughes, G.~Scovazzi, and L.P. Franca.
\newblock {\em {M}ultiscale and stabilized methods}.
\newblock Wiley, 2004.

\bibitem{john:2010}
V.~John and A.~Kindl.
\newblock Numerical studies of finite element {V}ariational {M}ultiscale
  {M}ethods for turbulent flow simulations.
\newblock {\em Computer Methods in Applied Mechanics and Engineering},
  199:841--852, 2010.

\bibitem{john:2007}
V.~John and M.~Roland.
\newblock Simulations of the turbulent channel flow at ${R}e_{\tau}=180$ with
  projection-based finite element {V}ariational {M}ultiscale {M}ethods.
\newblock {\em International Journal of Numerical Methods in Fluids},
  55:407--429, 2007.

\bibitem{knight:1998}
D.~Knight, G.~Zhou, N.~Okong'o, and V.Shukla.
\newblock Compressible {L}arge {E}ddy {S}imulation using unstructured grids.
\newblock Technical Report 98-0535, American Institute of Aeronautics and
  Astronautics, 1998.

\bibitem{koobus:2004}
B.~Koobus and C.~Farhat.
\newblock A {V}ariational {M}ultiscale method for the {L}arge {E}ddy
  {S}imulation of compressible turbulent flows on unstructured
  meshes----application to vortex shedding.
\newblock {\em Computer Methods in Applied Mechanics and Engineering},
  193:1367--1383, 2004.

\bibitem{lander:1997}
J.~Lander and B.J. Hoskins.
\newblock Believable scales and parameterizations in a spectral transform
  model.
\newblock {\em Monthly Weather Review}, 125:292--303, 1997.

\bibitem{landmann:2008}
B.~Landmann, M.~Kessler, S.~Wagner, and E.~Kr{\"a}mer.
\newblock A parallel, high-order discontinuous {G}alerkin code for laminar and
  turbulent flows.
\newblock {\em Computers \& Fluids}, 37:427--438, 2008.

\bibitem{lenormand:2000}
E.~Lenormand, P.Sagaut, and L.~{Ta Phuoc}.
\newblock Large {E}ddy {S}imulation of subsonic and supersonic channel flow at
  moderate {R}eynolds number.
\newblock {\em International Journal of Numerical Methods in Fluids},
  32:369--406, 2000.

\bibitem{maggioni:2012}
A.~Maggioni.
\newblock Formulazione {DG-LES} per flussi turbolenti comprimibili: modelli e
  validazione in un canale piano.
\newblock Master's thesis, School of Industrial Engineering, Politecnico di
  Milano, 2012.

\bibitem{piomelli:2000}
M.~Pino Martin, U.~Piomelli, and G.V. Candler.
\newblock {Subgrid-Scale Models for Compressible {L}arge-{E}ddy {S}imulations}.
\newblock {\em Theoretical and Computational Fluid Dynamics}, 13:361--376,
  2000.

\bibitem{germano:2014}
M.Germano, A.~Abb\`a, R.~Arina, and L.~Bonaventura.
\newblock On the extension of the eddy viscosity model to compressible flows.
\newblock {\em Physics of Fluids}, 26(4):041702, 2014.

\bibitem{moser:1999}
R.D. Moser, J.~Kim, and N.N. Mansour.
\newblock Direct numerical simulation of turbulent channel flow up to
  $re_\tau=590$.
\newblock {\em Physics of Fluids}, 11:943--945, 1999.

\bibitem{munts:2007}
E.A. Munts, S.J. Hulshoff, and R.~de~Borst.
\newblock A modal-based multiscale method for large eddy simulation.
\newblock {\em Journal of Computational Physics}, 224:389--402, 2007.

\bibitem{restelli:2009}
M.~Restelli and F.X. Giraldo.
\newblock A conservative {D}iscontinuous {G}alerkin semi-implicit formulation
  for the {N}avier-{S}tokes equations in nonhydrostatic mesoscale modeling.
\newblock {\em SIAM Journal of Scientific Computing}, 31:2231--2257, 2009.

\bibitem{sagaut:2006}
P.~Sagaut.
\newblock {\em Large {E}ddy {S}imulation for {I}ncompressible {F}lows: {A}n
  {I}ntroduction.}
\newblock Springer Verlag, 2006.

\bibitem{schlichting:1979}
H.~Schlichting.
\newblock {\em Boundary-layer theory.7th edition}.
\newblock McGraw-Hill, 1979.

\bibitem{schmitt:2007}
F.G. Schmitt.
\newblock About {B}oussinesq's turbulent viscosity hypothesis: historical
  remarks and a direct evaluation of its validity.
\newblock {\em Comptes Rendus M{\'e}canique}, 335:617--627, 2007.

\bibitem{schulze:2009}
J.~C. Schulze, P.~J. Schmid, and J.~L. Sesterhenn.
\newblock Exponential time integration using {K}rylov subspaces.
\newblock {\em International Journal of Numerical Methods in Fluids},
  60:591--609, 2009.

\bibitem{lilly:1993}
A.~Scotti, C.~Meneveau, and D.~Lilly.
\newblock {Generalized {S}magorinsky Model for Anisotropic Grids}.
\newblock {\em Physics of Fluids}, 5(9):2306--2308, 1993.

\bibitem{sengupta:2007}
K.~Sengupta, F.~Mashayek, and G.B. Jacobs.
\newblock Large {E}ddy {S}imulation using a discontinuos {G}alerkin spectral
  method.
\newblock In {\em 45th AIAA Aerospace Sciences Meeting and Exhibit}. AIAA,
  AIAA-2007-402 2007.

\bibitem{spiteri:2002}
R.J. Spiteri and S.J. Ruuth.
\newblock {A New Class of Optimal High-Order {S}trong-{S}tability-{P}reserving
  Time Discretization Methods}.
\newblock {\em SIAM Journal of Numerical Analysis}, 40:469--491, 2002.

\bibitem{tumolo:2013}
G.~Tumolo, L.~Bonaventura, and M.~Restelli.
\newblock A semi-implicit, semi-{L}agrangian, p-adaptive {D}iscontinuous
  {G}alerkin method for the shallow water equations.
\newblock {\em Journal of Computational Physics}, 232:46Ð67, 2013.

\bibitem{uranga:2011}
A.~Uranga, P.O. Persson, M.~Drela, and J.~Peraire.
\newblock Implicit {L}arge {E}ddy {S}imulation of transition to turbulence at
  low {R}eynolds numbers using a {D}iscontinuous {G}alerkin method.
\newblock {\em International Journal for Numerical Methods in Engineering},
  87:232--261, 2011.

\bibitem{vanderbos:2010}
F.~van~der Bos and B.J. Geurts.
\newblock Computational error-analysis of a {D}iscontinuous {G}alerkin
  discretization applied to large-eddy simulation of homogeneous turbulence.
\newblock {\em Computer Methods in Applied Mechanics and Engineering},
  199:903--915, 2010.

\bibitem{vreman:1995}
B.~Vreman, B.J. Geurts, and H.~Kuerten.
\newblock .subgrid-modeling in {LES} of compressible flow.
\newblock {\em Applied Scientific Research}, 54:191--203, 1995.

\bibitem{wei:2011}
L.~Wei and A.~Pollard.
\newblock Direct numerical simulation of compressible turbulent channel flows
  using the {D}iscontinuous {G}alerkin method.
\newblock {\em Computers and Fluids}, 47:85--100, 2011.

\bibitem{yang:1993}
K.S. Yang and J.H. Ferziger.
\newblock Large-{E}ddy {S}imulation of turbulent obstacle flow using a dynamic
  subgrid-scale model.
\newblock {\em AIAA Journal}, 31:1406--1413, 1993.

\bibitem{yoshizawa:1986}
A.~Yoshizawa.
\newblock Statistical theory for compressible turbulent flows with the
  application to subgrid modeling.
\newblock {\em Physics of Fluids}, 29:2152--2164, 1986.

\bibitem{zang:1993}
Y.~Zang, R.L. Street, and J.R. Koseff.
\newblock A dynamic mixed subgrid-scale model and its application to turbulent
  recirculating flows.
\newblock {\em Physics of Fluids}, 5:3186--3196, 1993.

\end{thebibliography}
\end{document}